\documentclass[12pt]{article}
\usepackage{amssymb}
\usepackage{amsmath}
\usepackage{amsthm}
\usepackage{color}
\usepackage{comment}
\usepackage{graphicx}
\begin{document}

\def\cA{\mathcal{A}}
\def\cB{\mathcal{B}}
\def\bx{{\bf x}}
\def\bW{{\bf W}}
\def\ee{\varepsilon}

\newcommand{\removableFootnote}[1]{}

\newtheorem{theorem}{Theorem}
\newtheorem{lemma}[theorem]{Lemma}

\title{
On resolving singularities of piecewise-smooth discontinuous vector fields via small perturbations.
}
\author{
D.J.W.~Simpson\\
Institute of Fundamental Sciences\\
Massey University\\
Palmerston North\\
New Zealand}
\maketitle




\begin{abstract}

A two-fold singularity is a point on a discontinuity surface of a piecewise-smooth vector field
at which the vector field is tangent to the surface on both sides.
Due to the double tangency, forward evolution from a two-fold is typically ambiguous.
This is an especially serious issue for two-folds
that are reached by the forward orbits of a non-zero measure set of initial points.
The purpose of this paper is to explore the concept of
perturbing the vector field so that forward evolution is well-defined,
and characterising the perturbed dynamics in the limit that the size of the perturbation tends to zero.
This concept is applied to a two-fold in two dimensions.
Three forms of perturbation: hysteresis, time-delay, and noise, are analysed individually.
In each case, the limit leads to a novel probabilistic notion of forward evolution from the two-fold.

\end{abstract}

\section{Introduction}
\label{sec:INTRO}
\setcounter{equation}{0}

\begin{figure}[t!]
\begin{center}
\setlength{\unitlength}{1cm}
\begin{picture}(11.2,6.6)
\put(0,3.5){\includegraphics[height=2.7cm]{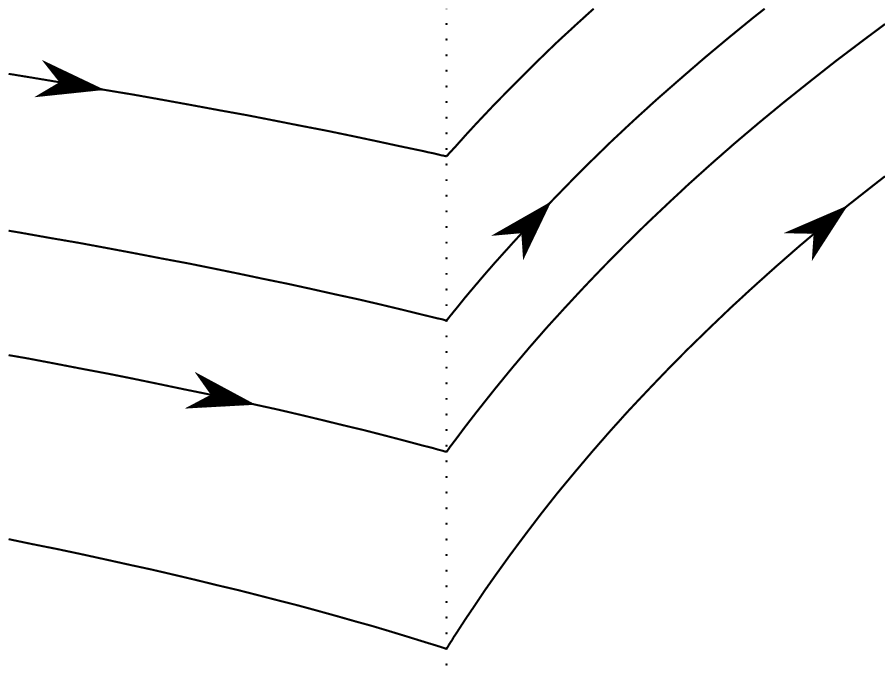}}	
\put(3.8,3.5){\includegraphics[height=2.7cm]{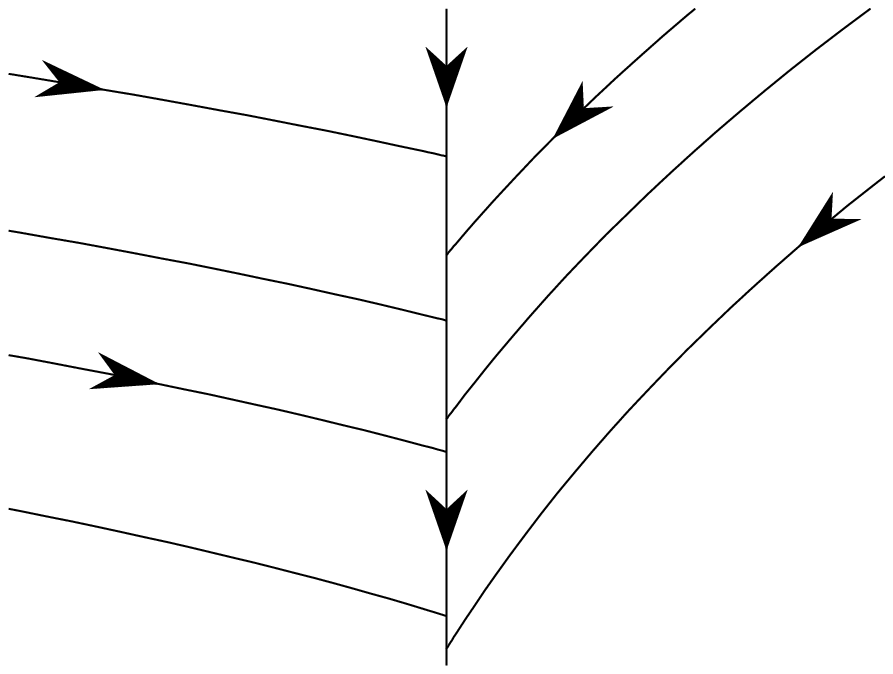}}
\put(7.6,3.5){\includegraphics[height=2.7cm]{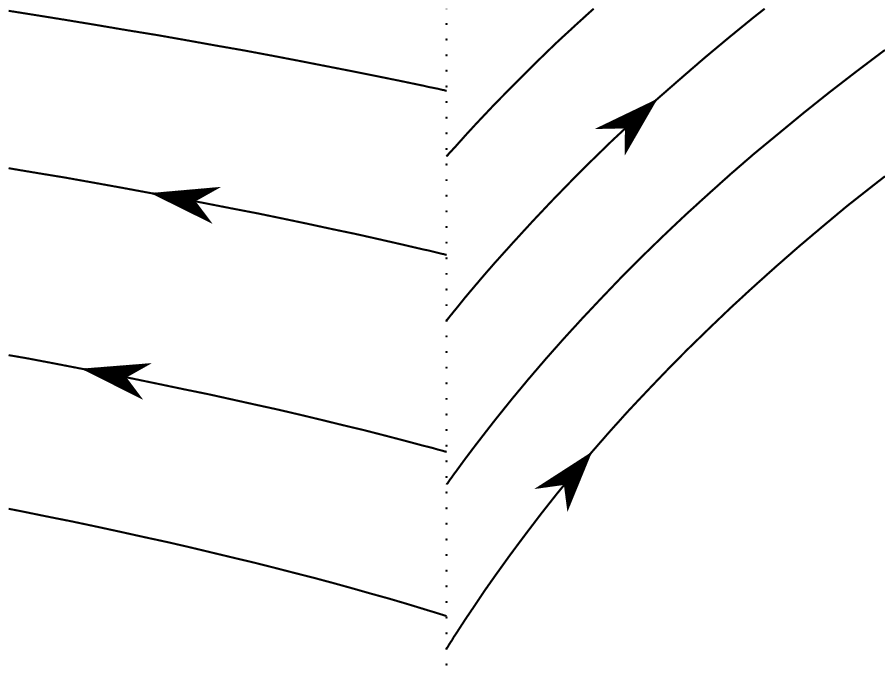}}
\put(1.7,0){\includegraphics[height=2.7cm]{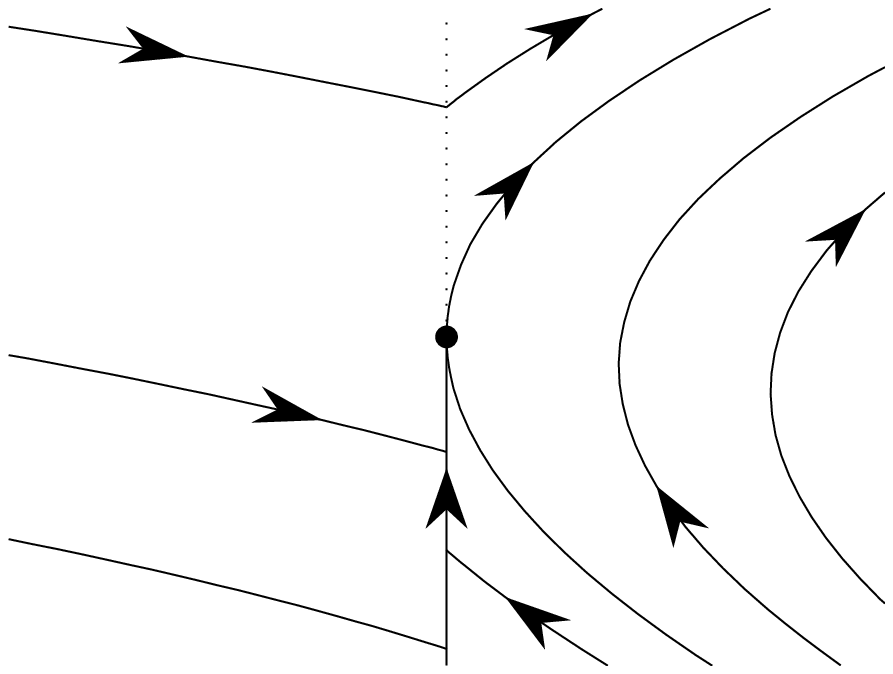}}
\put(5.9,0){\includegraphics[height=2.7cm]{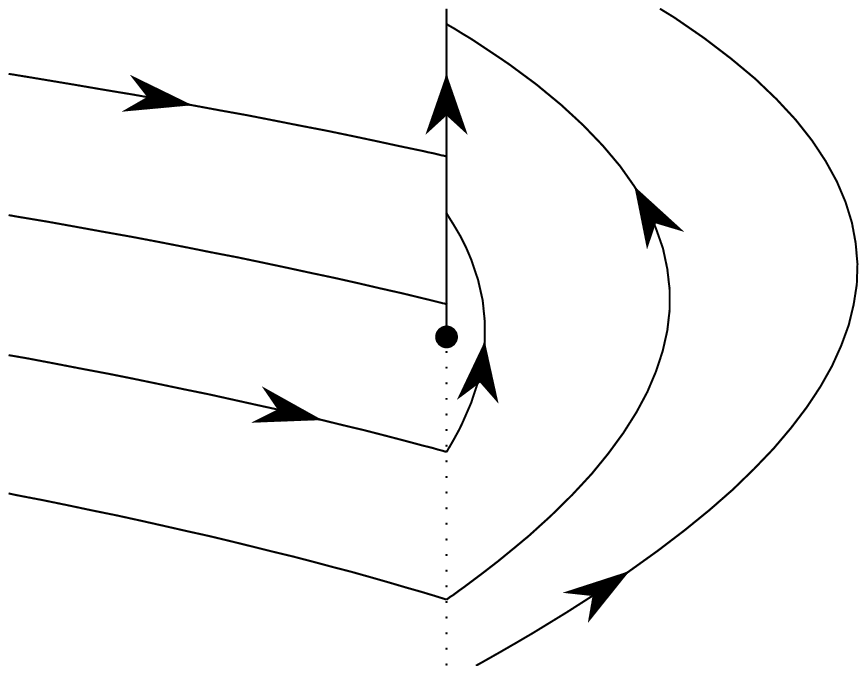}}
\put(1.1,6.5){\small crossing}
\put(4.5,6.5){\small stable sliding}
\put(8.6,6.5){\small repelling}
\put(2.2,3){\small visible tangency}
\put(6.1,3){\small invisible tangency}
\end{picture}
\caption{
Schematic phase portraits of a system of autonomous,
piecewise-smooth, discontinuous ODEs near a discontinuity surface.
The top sketches illustrate the three cases for the dynamics local
to a point on a discontinuity surface that is generic.
The bottom sketches show tangencies which are examples of non-generic points on a discontinuity surface.
\label{fig:schemFilippov}
}
\end{center}
\end{figure}

Systems of piecewise-smooth, discontinuous ODEs are utilised in many areas to model phenomena involving switching events.
Except in special cases, the dynamical behaviour local to a point on a discontinuity surface
conforms to one of three scenarios: crossing, stable sliding, or repelling.
These are illustrated in Fig.~\ref{fig:schemFilippov}.
In the first scenario trajectories simply cross the discontinuity surface.
In the second scenario the discontinuity surface attracts nearby trajectories.
Here Filippov's solution \cite{Fi60,Fi88,DiBu08,LeNi04} may be used
to define {\em sliding motion} on the discontinuity surface.
Sliding motion has been described in models of diverse phenomena, including
oscillators subject to dry friction \cite{BlCz99,WiDe00,OeHi96,FeMo94},
relay control \cite{Jo03,JoRa99,DiJo01},
and population dynamics \cite{DeGr07,DeFe06,AmOl12,TaLi12}.
In the third scenario trajectories head away from the discontinuity surface from both sides.
In this case forward evolution from a point on the discontinuity surface is ambiguous or non-unique
and a differential inclusion may be used to define a set-valued solution \cite{De92,Sm00,Co08c}.
However, from an applied viewpoint this is typically not an important issue
because there are no initial conditions for which forward evolution leads to the discontinuity surface.

Points on discontinuity surfaces at which the vector field is tangent to the surface on one side
usually represent boundaries at which the nature of the discontinuity surface changes
between one of the three generic scenarios.
For example, at the visible tangency shown in Fig.~\ref{fig:schemFilippov},
forward evolution changes from sliding motion to regular motion along the tangent trajectory.
Here the tangency is referred to as {\em visible} because the tangent trajectory is a solution of the system.
In contrast, the tangency shown in the bottom right sketch of Fig.~\ref{fig:schemFilippov} is said to be {\em invisible}
because the tangent trajectory is not seen.

\begin{figure}[b!]
\begin{center}
\setlength{\unitlength}{1cm}
\begin{picture}(8,6)
\put(0,0){\includegraphics[height=6cm]{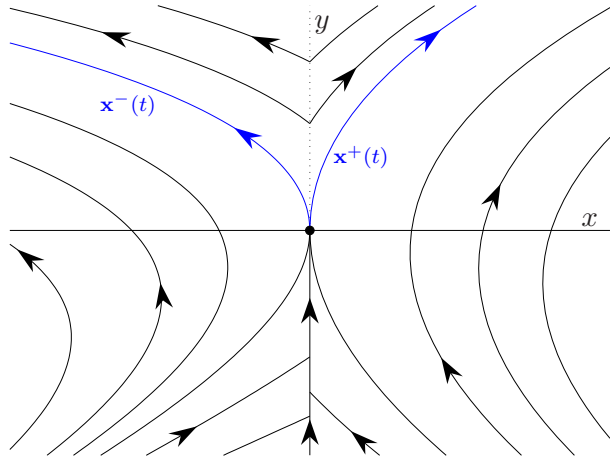}}
\put(7.6,3.07){\small $x$}
\put(4.05,5.7){\small $y$}
\put(1.2,4.6){\scriptsize \color{blue} $\bx^-(t)$}
\put(4.3,3.9){\scriptsize \color{blue} $\bx^+(t)$}
\end{picture}
\caption{
Dynamics local to a generic planar two-fold for which both tangent trajectories,
$\bx^+(t)$ and $\bx^-(t)$, are visible and point in the same direction at the two-fold.
\label{fig:schemPlanarTwoFold}
}
\end{center}
\end{figure}

Points on discontinuity surfaces at which the vector field is tangent to the surface on both sides are known as {\em two-folds}.
Such points arise generically in systems of at least three dimensions \cite{Te90,Te93},
or as codimension-one bifurcations in planar systems \cite{KuRi03}.
As with repelling discontinuity surfaces, forward evolution from two-folds is generally ambiguous.
However, unlike for repelling discontinuity surfaces,
the ambiguity of forward evolution from a two-fold 
is a serious issue if many trajectories go into the two-fold.
This is the case for the two-fold depicted in Fig.~\ref{fig:schemPlanarTwoFold},
because the area of the set of all points whose forward orbits intersect the two-fold is non-zero.

For planar systems there are several distinct generic two-folds.
These may be distinguished by whether neither, one, or both 
tangent trajectories are visible,
and the relative direction of these trajectories at the two-fold \cite{KuRi03,GuSe11}\removableFootnote{
Thus there are six distinct generic two-folds in two dimensions.
These each have a unique unfolding, except the anti-collinear visible-invisible two-fold which has two distinct unfoldings.
I reckon that I don't need to say this because I haven't said
exactly how I'm ``distinguishing'' different two-folds.
Unfoldings of codimension-three scenarios in planar Filippov systems are given in \cite{BuDe12}.
}.
For Fig.~\ref{fig:schemPlanarTwoFold}, both tangent trajectories are visible
and they point in the same direction at the two-fold.
The tangent trajectories are the orbits of the left and right vector fields that pass through the origin.
We write them as, $\bx^\pm(t) \equiv \left( x^\pm(t), y^\pm(t) \right)$, and for simplicity assume, $\bx^\pm(0) = (0,0)$.

In three dimensions, a two-fold for which at least one tangent trajectory is visible,
may exhibit an ambiguity similar to that of Fig.~\ref{fig:schemPlanarTwoFold}
in that the forward orbits of a nonzero volume of initial points intersect the two-fold \cite{Te90,Te93}.
Dynamics local to generic two-folds in three dimensions have been systematically classified \cite{JeCo09,CoJe11,FeGa12}.
In more than three dimensions, two-folds have similar properties \cite{CoJe12}.
Also, it has recently been shown that two-folds arise in models of simple circuit systems \cite{DiCo11},
and have deep connections to {\em folded nodes} of slow-fast systems \cite{DeJe11,DeJe12}.

However, dynamical behaviour near two-folds should only be understood with the knowledge that
as mathematical models nonsmooth equations only represent an approximation to reality.
Discontinuities are used to model events, such as impacts,
that on an appropriately fine scale are not instantaneous and perhaps actually smooth, although highly nonlinear.
Model inaccuracies are crucial at sensitive regions of phase space;
indeed two-folds are points of extreme sensitivity.

To address this issue, in this paper physically motivated perturbations of a discontinuous ODE system
are analysed with the aim of resolving the ambiguity of forward evolution at two-folds.
In order to obtain analytical results the scenario depicted in Fig.~\ref{fig:schemPlanarTwoFold} is studied
because it is the simplest two-fold that exhibits the serious ambiguity issue of interest.
Higher-dimensional scenarios are left for future work.

The following three types of perturbation are considered separately\removableFootnote{
A smoothing of the discontinuity is a fourth type of perturbation.
Here the stable sliding region is perturbed into a stable slow manifold
and the unstable sliding region is perturbed into an unstable slow manifold.
For small $\ee > 0$, the two manifolds almost coincide.
Whether an orbit heads left or right is determined by which side of the manifold the orbit starts from.
It is not clear what is gained by this description.
The three perturbations described here give more interesting results because
each perturbation creates an intertwining of orbits about the perturbed stable sliding region.
}.
{\em Hysteresis} -- a common feature of control systems \cite{Ts84} --
is incorporated by replacing the discontinuity surface with a hysteretic band of width $2 \ee$.
In this context, and throughout the paper, $\ee > 0$ is assumed to be small.
Second, {\em time-delay} is added.
Specifically, it is supposed that the functional form of the equations governing forward evolution switches
at a time $\ee$ after trajectories intersect the discontinuity surface.
Small time-delay is inherent in many switching systems, such as relay control systems \cite{Ts84},
and represents the time lag between when the switching threshold is crossed and a change in dynamics is enacted.
Lastly, {\em noise} of amplitude $\ee$ is added to the differential equations.
Noise and uncertainty are ubiquitous in real-world systems.
Additive white Gaussian noise is used here as this is possibly
the simplest manner by which noise may be incorporated into the system.

With a small hysteretic band, stable sliding motion approaching the two-fold is replaced by rapid switching motion.
Importantly, for any $\ee > 0$, forward evolution is uniquely determined in a neighbourhood of the two-fold.
Forward orbits pass close to the two-fold and then are directed away from the two-fold
along a path close to either $\bx^+(t)$ or $\bx^-(t)$.
We are not interested in the precise path that orbits take,
only whether they eventually ``head right'',
that is follow close to $\bx^+(t)$ after passing by the two-fold, or ``head left''.
Given, say, an interval of initial points on the discontinuity surface below the two-fold,
for any $\ee > 0$, some fraction of forward orbits of these points head right after passing by the two-fold.
In the limit, $\ee \to 0$, this fraction approaches a constant value, $\varrho_{\rm hy}$.
When the system is instead perturbed by the inclusion of time-delay,
the dynamics fit a similar description,
but in the limit, $\ee \to 0$, the fraction of forward orbits that head right approaches
a different value, $\varrho_{\rm td}$.
The addition of noise yields a system of stochastic differential equations
that have a unique probabilistic solution for any $\ee > 0$.
From a suitable initial point below the two-fold,
the stochastic system defines a probability that sample paths from this point eventually head right.
As $\ee \to 0$, this probability approaches a value, $\varrho_{\rm no}$.
Thus in each case, with $\ee > 0$ forward evolution is well-defined,
and upon passing near the two-fold most orbits (or sample paths in the case of the noise)
follow a path close to $\bx^+(t)$ or $\bx^-(t)$.
In each case we obtain a limiting fraction or probability, $\varrho$, of trajectories that head right.

The concept of adding noise to remove an ambiguity is not new.
Under general assumptions the addition of noise to a system of nonsmooth ODEs with multiple solutions
gives a stochastic system with a unique solution \cite{Fl11,KrRo05,Da07}.
In the zero-noise limit, this solution may represent a useful
weighted selection of all possible solutions \cite{FlLa08,BoSu10}.
Recently such ideas have been extended and applied to transport equations
for fluid dynamics \cite{Fl09}.
However, different ways of adding noise can give different solutions in the zero-noise limit \cite{AtFl09}\removableFootnote{
I don't understand the results of this paper.
}.
In addition, for a slow-fast system, when both the zero-noise limit and slow-fast limit are taken,
the limiting solution may depend on the order in which the limits are applied \cite{Sa83}.
A complete description of the zero-noise limit has been achieved for some simple one-dimensional systems \cite{Ve83,BaBa82}.

Variations on this concept have been applied in alternate contexts.
For instance, forward evolution from an unstable hyperbolic equilibrium of an $N$-dimensional system of ODEs is motionless,
but with arbitrarily small noise, sample paths escape the equilibrium and do so, with high probability,
along a path near to the invariant manifold associated with the largest eigenvalue of the equilibrium \cite{Ki81}.
For systems of heteroclinic networks, trajectories have a ``choice'' at saddle equilibria as to which heteroclinic path to take.
The addition of noise can be used to replace this choice with a probability \cite{Ba11,ArSt03}.

The remainder of the paper is organised as follows.
Equations describing the planar two-fold of Fig.~\ref{fig:schemPlanarTwoFold} are introduced in \S\ref{sec:TWOFOLD}.
Perturbation by hysteresis, time-delay and noise
are described and analysed in \S\ref{sec:HYST}, \S\ref{sec:DELAY}, and \S\ref{sec:NOISE}, respectively.
Conclusions are given in \S\ref{sec:CONC}.

\section{A planar visible-visible two-fold}
\label{sec:TWOFOLD}
\setcounter{equation}{0}

We consider the dynamics of a general two-dimensional system of ODEs near a smooth discontinuity surface.
By choosing coordinates, $(X,Y)$, such that the discontinuity surface
coincides with the $Y$-axis, a description of the local dynamics of the system may be written as
\begin{equation}
\left[ \begin{array}{c} \dot{X} \\ \dot{Y} \end{array} \right] =
\left\{ \begin{array}{lc}
\left[ \begin{array}{c} f^{(L)}(X,Y) \\ g^{(L)}(X,Y) \end{array} \right] \;, & X < 0 \\
\left[ \begin{array}{c} f^{(R)}(X,Y) \\ g^{(R)}(X,Y) \end{array} \right] \;, & X > 0
\end{array} \right. \;,
\label{eq:odeOriginal}
\end{equation}
where we assume that the four functions that appear in (\ref{eq:odeOriginal})
are smooth on the closure of their associated half-planes.
In order for (\ref{eq:odeOriginal}) to exhibit the
two-fold of Fig.~\ref{fig:schemPlanarTwoFold} at the origin, we must have
\begin{equation}
\begin{split}
f^{(L)}(0,0) = 0 \;, \qquad
f^{(R)}(0,0) = 0 \;, \\
\frac{\partial f^{(L)}}{\partial Y}(0,0) < 0 \;, \qquad
\frac{\partial f^{(R)}}{\partial Y}(0,0) > 0 \;, \\
g^{(L)}(0,0) > 0 \;, \qquad
g^{(R)}(0,0) > 0 \;.
\end{split}
\label{eq:planarTwoFoldConditions}
\end{equation}
The first two equations of (\ref{eq:planarTwoFoldConditions})
specify that the left and right vector fields are tangent to the discontinuity surface at the origin,
the second two equations specify that both tangent trajectories are visible,
and the third two equations specify that both tangent trajectories are directed upward.

In order to minimise the number of parameters that must be considered, we define the scaled state variables
\begin{equation}
x = \frac{1}{\frac{\partial f^{(R)}}{\partial Y}(0,0) g^{(R)}(0,0)} X \;, \qquad
y = \frac{1}{g^{(R)}(0,0)} Y \;.
\label{eq:xy}
\end{equation}
Since only dynamics near the origin is of interest,
we expand the two halves of the vector field as Taylor series centred about the origin,
using $O(k)$ to denote terms that are order $k$, or larger, in $x$ and $y$.
Specifically, (\ref{eq:odeOriginal}) may be rewritten as
\begin{equation}
\left[ \begin{array}{c} \dot{x} \\ \dot{y} \end{array} \right] =
\left\{ \begin{array}{lc}
\left[ \begin{array}{c} -\cA y + O(|x|) + O(2) \\ \cB + O(1) \end{array} \right] \;, & x < 0 \\
\left[ \begin{array}{c} y + O(|x|) + O(2) \\ 1 + O(1) \end{array} \right] \;, & x > 0
\end{array} \right. \;,
\label{eq:ode}
\end{equation}
where
\begin{equation}
\cA = -\frac{\frac{\partial f^{(L)}}{\partial Y}(0,0)}
{\frac{\partial f^{(R)}}{\partial Y}(0,0)} \;, \qquad
\cB = \frac{g^{(L)}(0,0)}{g^{(R)}(0,0)} \;.
\label{eq:AB}
\end{equation}
Consequently, $\cA$ and $\cB$ may take any positive values.

The limiting fractions or probabilities, $\varrho_{\rm hy}$, $\varrho_{\rm td}$ and $\varrho_{\rm no}$,
determined in the the following sections,
depend only on $\cA$ and $\cB$, and not on any higher order terms.
For this reason we write $\varrho = \varrho(\cA,\cB)$
for each type of perturbation.
Furthermore, in each case $\varrho$ satisfies the symmetry property
\begin{equation}
\varrho \left( \frac{1}{\cA},\frac{1}{\cB} \right) = 1 - \varrho(\cA,\cB) \;.
\label{eq:varrhoSym}
\end{equation}
To see this, note that if we were to scale $x$ and $y$ in (\ref{eq:ode}) such that the vector field for $x < 0$ is simply
$(\dot{x},\dot{y}) = (-y,1)$, plus higher order terms,
then the vector field for $x > 0$ would become
$(\dot{x},\dot{y}) = \left( \frac{1}{\cA} y, \frac{1}{\cB} \right)$, plus higher order terms.
Reversing the sign of $x$ would then recover (\ref{eq:ode}),
but with $\cA$ and $\cB$ replaced by $\frac{1}{\cA}$ and $\frac{1}{\cB}$,
and the concepts of left and right switched.
Therefore if $\varrho(\cA,\cB)$ represents the limiting fraction or probability of trajectories heading right,
then $\varrho \left( \frac{1}{\cA},\frac{1}{\cB} \right)$ is the probability of trajectories heading left.
Assuming these probabilities sum to $1$ (which is the case for each of the three types of perturbation),
we therefore have (\ref{eq:varrhoSym}).

\section{Hysteresis}
\label{sec:HYST}
\setcounter{equation}{0}

In this section the switching condition of the system (\ref{eq:ode}) is replaced
with a hysteretic band spanning $-\ee \le x \le \ee$.
More precisely, suppose that the equations governing forward evolution switch from
the right half-system to the left half-system when an orbit reaches $x = -\ee$, 
and similarly the governing equations switch from
the left half-system to the right half-system when an orbit reaches $x = \ee$.
Algebraically this equates to
\begin{equation}
\left[ \begin{array}{c} \dot{x} \\ \dot{y} \end{array} \right] =
\left\{ \begin{array}{lc}
\left[ \begin{array}{c} -\cA y + O(|x|) + O(2) \\ \cB + O(1) \end{array} \right] \;, & {\rm until~} x = \ee \\
\left[ \begin{array}{c} y + O(|x|) + O(2) \\ 1 + O(1) \end{array} \right] \;, & {\rm until~} x = -\ee
\end{array} \right. \;.
\label{eq:hysteresis}
\end{equation}
In \cite{DiJo02}, an analogous $O(\ee)$ hysteretic band was added to a
piecewise-smooth system exhibiting a stable periodic orbit with sliding motion.
Here the authors found that the periodic orbit persisted only for relatively small values of $\ee$.
More generally, hysteresis may be the cause of complicated and chaotic dynamics \cite{GoMe01,KaWa11,VaGe01}.

\begin{figure}[b!]
\begin{center}
\setlength{\unitlength}{1cm}
\begin{picture}(5.33,12)
\put(0,0){\includegraphics[height=12cm]{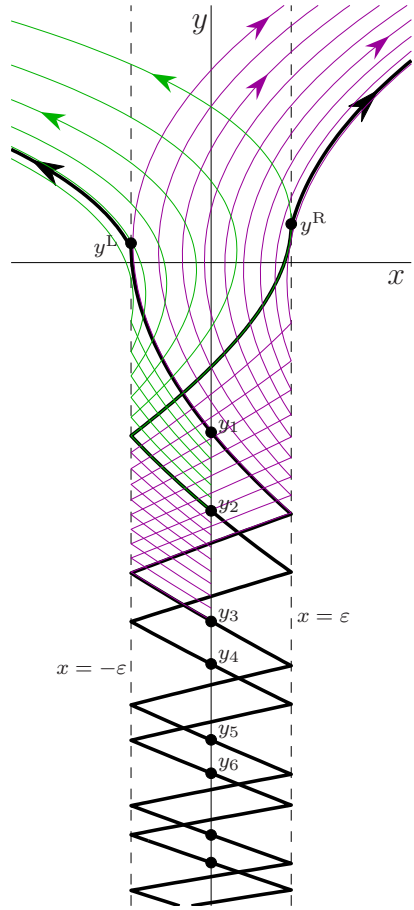}}
\put(5,8.3){$x$}
\put(2.4,11.7){$y$}
\put(3.84,9){\scriptsize $y^{\rm R}$}
\put(1.1,8.7){\scriptsize $y^{\rm L}$}
\put(2.75,6.36){\scriptsize $y_1$}
\put(2.75,5.3){\scriptsize $y_2$}
\put(2.75,3.85){\scriptsize $y_3$}
\put(2.75,3.3){\scriptsize $y_4$}
\put(2.75,2.3){\scriptsize $y_5$}
\put(2.75,1.86){\scriptsize $y_6$}
\put(.6,3.1){\scriptsize $x = -\ee$}
\put(3.8,3.8){\scriptsize $x = \ee$}
\end{picture}
\caption{
A sketch of phase space near the origin of the hysteretic system (\ref{eq:hysteresis}).
The thick curves are the critical switching trajectories
that divide orbits that eventually head right and orbits that eventually head left.
For instance, as indicated, the forward orbit of any point, $(0,y)$, with $y_3 < y \le y_2$,
that initially follows the right half-system, eventually heads right,
whereas if $y_2 < y \le y_1$, then the orbit eventually heads left.
\label{fig:schemHysteresis}
}
\end{center}
\end{figure}

The stable sliding motion of (\ref{eq:ode}) that occurs at $x = 0$ for $y < 0$,
becomes rapid switching motion across $-\ee \le x \le \ee$ for the system (\ref{eq:hysteresis}), see Fig.~\ref{fig:schemHysteresis}.
Upon reaching the vicinity of the origin, some orbits head right, and the others head left.
For each orbit, there is no ambiguity associated with its forward evolution.
The purpose of this section is to precisely formulate, and then compute, the fraction of orbits that head right.

There are two critical trajectories that separate orbits by whether
they eventually head right or eventually head left.
These trajectories have a tangential intersection with one of the hysteretic switching manifolds $(x = \pm \ee)$.
The intersections occur at some points, $(-\ee,y^{\rm L})$ and $(\ee,y^{\rm R})$, as shown in Fig.~\ref{fig:schemHysteresis}.
For our purposes it is sufficient to note that the tangency points satisfy
\begin{equation}
y^{\rm L} = O(\ee) \;, \qquad
y^{\rm R} = O(\ee) \;.
\label{eq:yminusyplus}
\end{equation}
We follow the critical trajectories backward from their tangential intersections
assuming that switching occurs at every intersection with $x = \pm \ee$, as shown in Fig.~\ref{fig:schemHysteresis}.
The values, $y_k$, for $k \ge 1$, are used to denote the successive intersections of the critical trajectories with $x = 0$
at which the trajectories are governed by the right half-system.
By using odd values of $k$ for the critical trajectory passing through $(-\ee,y^{\rm L})$,
and even values of $k$ for the critical trajectory passing through $(\ee,y^{\rm R})$,
we have $y_k < y_{k-1}$, for all $k \ge 2$.
Now consider the forward orbit of a point, $(0,y)$, with $y_k < y \le y_{k-1}$ for some $k \ge 2$,
that initially follows the right half-system.
Simply, if $k$ is odd the orbit eventually heads right,
and if $k$ is even the orbit eventually heads left.

If we treat the initial point, $(0,y)$, as fixed, and decrease $\ee$ to zero,
the future of the forward orbit changes between eventually heading right and eventually heading left,
and in the limit $\ee \to 0$ its future is undefined.
For this reason, we instead consider a large collection of initial points
and study the fraction of points for which forward evolution eventually heads right.
To make this precise, we first let $x^* > \ee$ be an $\ee$-independent value chosen
sufficiently small to ensure that global dynamics do not play a role.
We then specify that an orbit ``heads right'' if the orbit first exits the region $|x| < x^*$, at $x = x^*$.
We define the following indicator function
for the forward orbit of any point, $(0,y)$, that initially follows the right half-system:
\begin{equation}
\sigma_\ee(y) = \left\{ \begin{array}{lc}
1 \;, & {\rm if~the~orbit~heads~right} \\
0 \;, & {\rm otherwise}
\end{array} \right. \;.
\end{equation}
(The assumption that orbits initially follow the right half-system is taken without loss of generality
and the value of $\varrho_{\rm hy}$ is not affected by this choice.)
For any interval of $y$-values, $I$, we let
\begin{equation}
q_I = \frac{\int_I \sigma_\ee(y) \,dy}{\int_I \,dy} \;,
\label{eq:qIdef}
\end{equation}
represent the fraction of the interval $I$ for which forward evolution eventually heads right.
Already we have shown
\begin{equation}
q_{[y_k,y_{k-2}]} = \frac{y_{k-1} - y_k}{y_{k-2} - y_k} \;, {\rm ~for~} k = 3,5,7,\ldots \;.
\label{eq:qk}
\end{equation}
Our goal is to derive
\begin{equation}
\varrho_{\rm hy} \equiv \lim_{\ee \to 0} q_{[y_{\rm min},y_{\rm max}]} \;,
\label{eq:varrhohysteresis}
\end{equation}
where $[y_{\rm min},y_{\rm max}]$ is a small $\ee$-independent interval of the $y$-axis below the origin.
In the following paragraphs we will see that $\varrho_{\rm hy}$ is independent of $y_{\rm min}$ and $y_{\rm max}$,
and thus $\varrho_{\rm hy}$ is an fundamental quantity of the unperturbed system.
Before we consider such $\ee$-independent intervals,
it is useful to first compute an explicit formula for $q_{[y_k,y_{k-2}]}$.

For any initial point $(x_0,y_0)$, with $x_0 = O(\ee)$ and $y_0 = O \left( \sqrt{\ee} \right)$,
if evolution (forward or backward) is governed purely by the left half-system up to a time $t = O \left( \sqrt{\ee} \right)$,
then the location of the orbit at this time is given by
\begin{equation}
\begin{split}
x^{(L)}(t) &= x_0 - \cA y_0 t - \frac{\cA \cB}{2} t^2 + O \left( \ee^{\frac{3}{2}} \right) \;, \\
y^{(L)}(t) &= y_0 + \cB t + O(\ee) \;.
\end{split}
\label{eq:xtytL}
\end{equation}
Similarly, if evolution is purely governed by the right half-system, then the orbit is located at
\begin{equation}
\begin{split}
x^{(R)}(t) &= x_0 + y_0 t + \frac{1}{2} t^2 + O \left( \ee^{\frac{3}{2}} \right) \;, \\
y^{(R)}(t) &= y_0 + t + O(\ee) \;.
\end{split}
\label{eq:xtytR}
\end{equation}
By substituting $(x_0,y_0) = (-\ee,y^{\rm L})$,
and $\left( x^{(R)}(t),y^{(R)}(t) \right) = (0,y_1)$ into (\ref{eq:xtytR}), we obtain
\begin{equation}
y_1 = -\sqrt{2} \sqrt{\ee} + O(\ee) \;.
\label{eq:y1}
\end{equation}
Similar substitutions using both (\ref{eq:xtytL}) and (\ref{eq:xtytR}) lead to
\begin{equation}
y_2 = -\sqrt{2} \sqrt{1 + \frac{2 \cB}{\cA}} \sqrt{\ee} + O(\ee) \;.
\label{eq:y2}
\end{equation}
In this fashion we obtain the recursion relation
\begin{equation}
y_k = -\sqrt{y_{k-2}^2 + 4 \left( 1 + \frac{\cB}{\cA} \right) \ee} + O(\ee) \;, {\rm ~for~} k \ge 3 \;.
\label{eq:ykplus2}
\end{equation}
Then it is easy to show that the combination
(\ref{eq:y1})-(\ref{eq:ykplus2}), gives the explicit formula
\begin{equation}
y_k = \left\{ \begin{array}{lc}
-\sqrt{2} \sqrt{k + \frac{(k-1) \cB}{\cA}} \sqrt{\ee} + O(\ee) \;, & {\rm ~for~} k = 1,3,5,\ldots \\
-\sqrt{2} \sqrt{k-1 + \frac{k \cB}{\cA}} \sqrt{\ee} + O(\ee) \;, & {\rm ~for~} k = 2,4,6,\ldots
\end{array} \right. \;.
\label{eq:yk}
\end{equation}
By substituting (\ref{eq:yk}) into (\ref{eq:qk}), and considering large $k$, we obtain
\begin{equation}
q_{[y_k,y_{k-2}]} = \frac{\cA}{\cA+\cB} + O \left( \frac{1}{k} \right) + O(\sqrt{\ee}) \;, {\rm ~for~} k = 3,5,7,\ldots \;.
\label{eq:qk2}
\end{equation}

\begin{figure}[t!]
\begin{center}
\setlength{\unitlength}{1cm}
\begin{picture}(9.73,7)
\put(.4,0){\includegraphics[height=7cm]{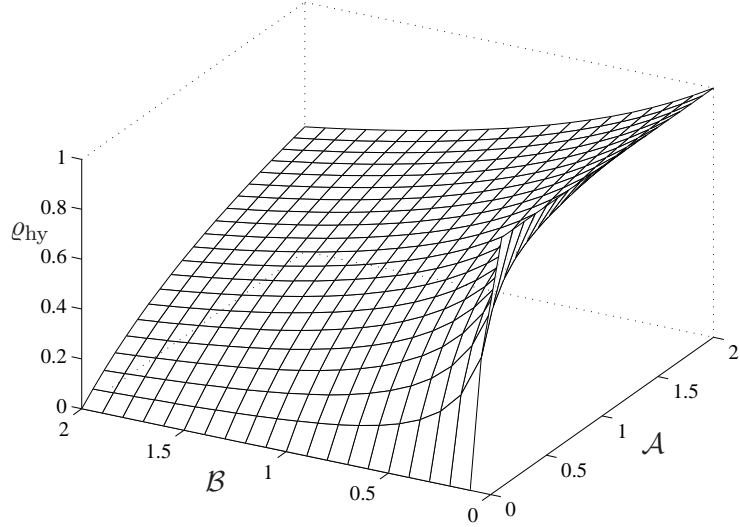}}
\put(8.4,1){\small $\cA$}
\put(2.6,.5){\small $\cB$}
\put(0,3.9){\small $\varrho_{\rm hy}$}
\end{picture}
\caption{
The fraction of orbits that head right after passing near the two-fold
for the hysteretic system (\ref{eq:hysteresis}) in the limit $\ee \to 0$,
$\varrho_{\rm hy}(\cA,\cB) = \frac{\cA}{\cA + \cB}$, (\ref{eq:varrhohysteresis2}).
\label{fig:ABhysteresis}
}
\end{center}
\end{figure}

Naturally we would like to use (\ref{eq:qk2}) to evaluate (\ref{eq:varrhohysteresis}).
Let $k_{\rm min}$ and $k_{\rm max}$ be the smallest and largest odd integers for which $y_k \in [y_{\rm min},y_{\rm max}]$.
Then, by (\ref{eq:qIdef}),
\begin{equation}
q_{[y_{\rm min},y_{\rm max}]} =
\frac{\sum_{k = k_{\rm min} + 2, k_{\rm min} + 4, \ldots}^{k_{\rm max}} (y_{k-1} - y_k)}
{y_{k_{\rm min}} - y_{k_{\rm max}}}
+ O(\ee) \;,
\label{eq:qI}
\end{equation}
where the error term corresponds to the $O(\ee)$ differences,
$y_{\rm max} - y_{k_{\rm min}}$ and $y_{k_{\rm max}} - y_{\rm min}$.
The quotient, (\ref{eq:qI}), may be written as a convex combination of,
$q_{[y_k,y_{k-2}]} = \frac{y_{k-1} - y_k}{y_{k-2} - y_k}$,
for all odd $k$ with $k_{\rm min}+2 \le k \le k_{\rm max}$.
For this reason we have the bound\removableFootnote{
Alternatively, this bound is an immediate consequence of the following result
\begin{lemma}
Given constants, $a_k, b_k > 0$, for $k = 1,2,\ldots$, for any $n \ge 1$,
\begin{equation}
{\rm min} \left[ \frac{a_1}{b_1} ,\ldots, \frac{a_n}{b_n} \right] \le
\frac{\sum_{k=1}^n a_k}{\sum_{k=1}^n b_k} \le
{\rm max} \left[ \frac{a_1}{b_1} ,\ldots, \frac{a_n}{b_n} \right] \;.
\end{equation}
\end{lemma}
For $n = 1$ the result is trivial.
For $n = 2$ the result may be proved via direct calculations.
For any $n \ge 3$ the result may be proved by induction.
}:
\begin{equation}
{\rm min}[q_{[y_k,y_{k-2}]}] \le
q_{[y_{\rm min},y_{\rm max}]} + O(\ee) \le
{\rm max}[q_{[y_k,y_{k-2}]}] \;,
\label{eq:qIbound}
\end{equation}
where the extrema are computed over all odd $k$ with $k_{\rm min}+2 \le k \le k_{\rm max}$.
However, $k_{\rm min}, k_{\rm max} \to \infty$ as $\ee \to 0$,
and since it has not been shown that the $O(\sqrt{\ee})$ term in (\ref{eq:qk2}) remains bounded as $k \to \infty$,
equations (\ref{eq:qk2}) and (\ref{eq:qIbound}) are insufficient for us to
derive the anticipated and intuitive result, (\ref{eq:qI2}), given below.
We have overcome this technical difficulty via an involved procedure of carefully bounding error terms.
We state the following result and provide a proof in Appendix \ref{sec:HYSTPROOF}.
\begin{lemma}
Given any $-1 \ll y_{\rm min} < y_{\rm max} < 0$,
\begin{equation}
q_{[y_{\rm min},y_{\rm max}]} \to \frac{\cA}{\cA+\cB} \;, {\rm ~as~} \ee \to 0 \;.
\label{eq:qI2}
\end{equation}
\label{le:qI}
\end{lemma}
The condition, $-1 \ll y_{\rm min}$, ensures that global features of the system do not influence the dynamics.
By (\ref{eq:varrhohysteresis}) and (\ref{eq:qI2}),
\begin{equation}
\varrho_{\rm hy} = \varrho_{\rm hy}(\cA,\cB) = \frac{\cA}{\cA + \cB} \;.
\label{eq:varrhohysteresis2}
\end{equation}
As mentioned above, $\varrho_{\rm hy}$ is independent of the values of $y_{\rm min}$ and $y_{\rm max}$.
A plot of $\varrho_{\rm hy}(\cA,\cB)$ is shown in Fig.~\ref{fig:ABhysteresis}\removableFootnote{
MJeffrey comments that he would like to see more numerics
for the hysteretic and time-delayed cases
but I don't see how this would add value.
}.

\section{Time-delay}
\label{sec:DELAY}
\setcounter{equation}{0}

Now suppose that the switching condition of (\ref{eq:ode}) involves a time-delay of length $\ee$:
\begin{equation}
\left[ \begin{array}{c} \dot{x} \\ \dot{y} \end{array} \right] =
\left\{ \begin{array}{lc}
\left[ \begin{array}{c} -\cA y + O(|x|) + O(2) \\ \cB + O(1) \end{array} \right] \;, & x(t-\ee) < 0 \\
\left[ \begin{array}{c} y + O(|x|) + O(2) \\ 1 + O(1) \end{array} \right] \;, & x(t-\ee) > 0
\end{array} \right. \;.
\label{eq:timedelay}
\end{equation}
The inclusion of time-delay may destroy periodic orbits with sliding \cite{DiJo02}.
Dynamics and bifurcations specific to 
general piecewise-smooth ODE systems with time-delayed switching conditions
are discussed in \cite{Si06}, and in \cite{SiKo10b} for systems with symmetry.
If both hysteresis and time-delay are incorporated in the switching condition, 
there are yet more potential complications \cite{CoDi07}.

\begin{figure}[t!]
\begin{center}
\setlength{\unitlength}{1cm}
\begin{picture}(5.33,12)
\put(0,0){\includegraphics[height=12cm]{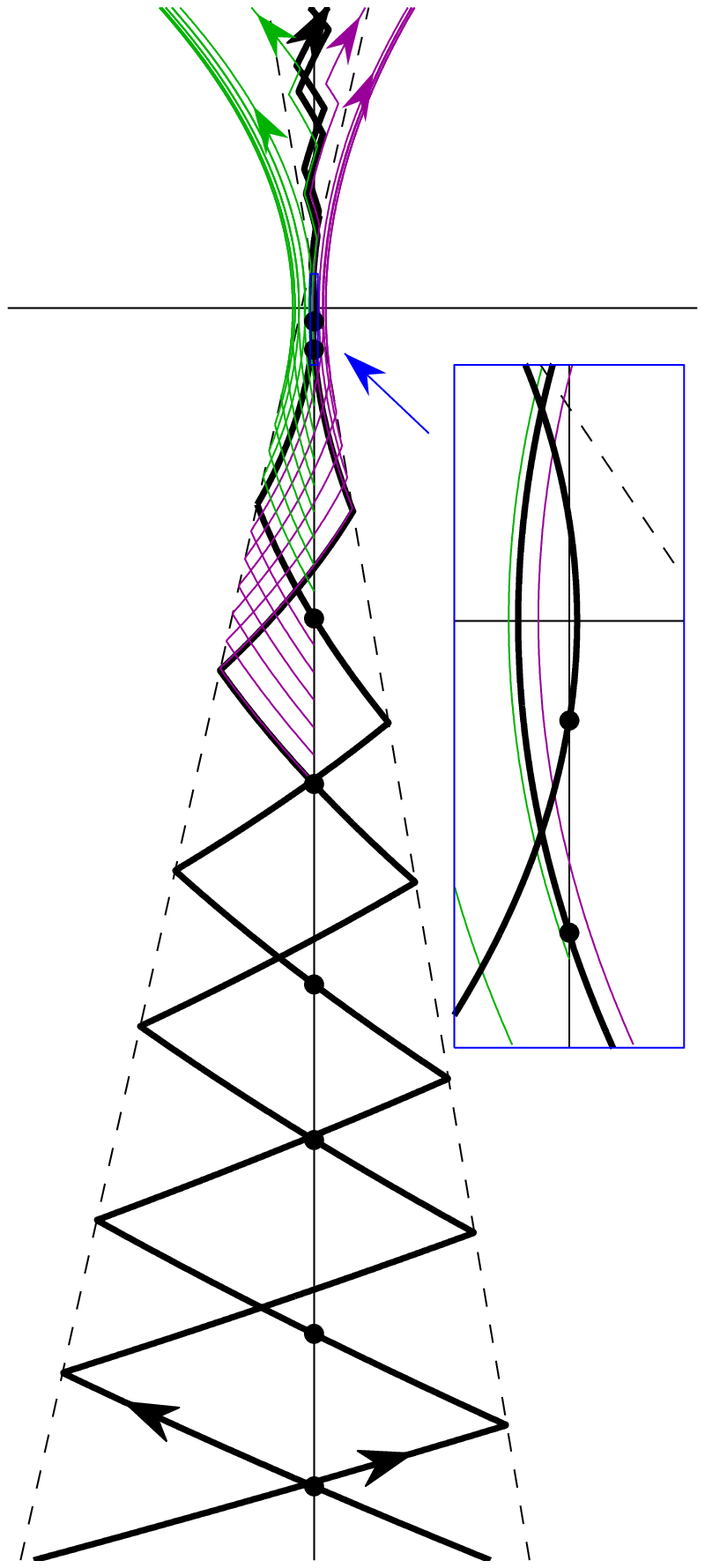}}
\put(5,9.75){$x$}
\put(0,1.7){$\Gamma^-$}
\put(3.8,1.7){$\Gamma^+$}
\put(4.92,7.02){\scriptsize $x$}
\put(4.4,8.8){\scriptsize $y$}
\put(4.5,4.9){\scriptsize $y_1$}
\put(4.45,6.5){\scriptsize $\tilde{y}_1$}
\put(2.45,7.35){\scriptsize $y_2$}
\put(2.53,5.95){\scriptsize $y_3$}
\put(2.47,4.5){\scriptsize $y_4$}
\put(2.4,2.85){\scriptsize $y_5$}
\put(2.46,1.82){\scriptsize $y_6$}
\put(2.4,.25){\scriptsize $y_7$}
\end{picture}
\caption{
\label{fig:schemTimeDelayLower}
A sketch of typical dynamics near the origin of the time-delayed system (\ref{eq:timedelay}).
}
\end{center}
\end{figure}

About the negative $y$-axis, dynamics of the time-delayed system (\ref{eq:timedelay}) 
switches at curves, $\Gamma^-$ and $\Gamma^+$
(shown in Fig.~\ref{fig:schemTimeDelayLower})
that may be thought of as the result of evolving the $y$-axis by a time $\ee$ under each half-system.
As in the previous section we use series solutions to obtain analytical results.
If $x_0 = O(\ee^2)$ and $y_0 = O(\ee)$,
then evolution up to a time, $t = O(\ee)$, for the left half-system of (\ref{eq:timedelay}) is given by
\begin{equation}
\begin{split}
x^{(L)}(t) &= x_0 - \cA y_0 t - \frac{\cA \cB}{2} t^2 + O(\ee^3) \;, \\
y^{(L)}(t) &= y_0 + \cB t + O(\ee^2) \;,
\end{split}
\label{eq:xtytL2}
\end{equation}
and for the right half-system of (\ref{eq:timedelay}), is given by
\begin{equation}
\begin{split}
x^{(R)}(t) &= x_0 + y_0 t + \frac{1}{2} t^2 + O(\ee^3) \;, \\
y^{(R)}(t) &= y_0 + t + O(\ee^2) \;.
\end{split}
\label{eq:xtytR2}
\end{equation}
From (\ref{eq:xtytL2}) and (\ref{eq:xtytR2}) we obtain
\begin{eqnarray}
\Gamma^- = \left\{ (x,y) ~\bigg|~ x = \ee y - \frac{1}{2} \ee^2 + O(\ee^3) \;,~y = O(\ee) \right\} \;,
\label{eq:Gammaminus} \\
\Gamma^+ = \left\{ (x,y) ~\bigg|~ x = -\ee \cA y + \frac{\cA \cB}{2} \ee^2 + O(\ee^3) \;,~y = O(\ee) \right\} \;.
\label{eq:Gammaplus}
\end{eqnarray}

\begin{figure}[t!]
\begin{center}
\setlength{\unitlength}{1cm}
\begin{picture}(8,12)
\put(0,0){\includegraphics[height=12cm]{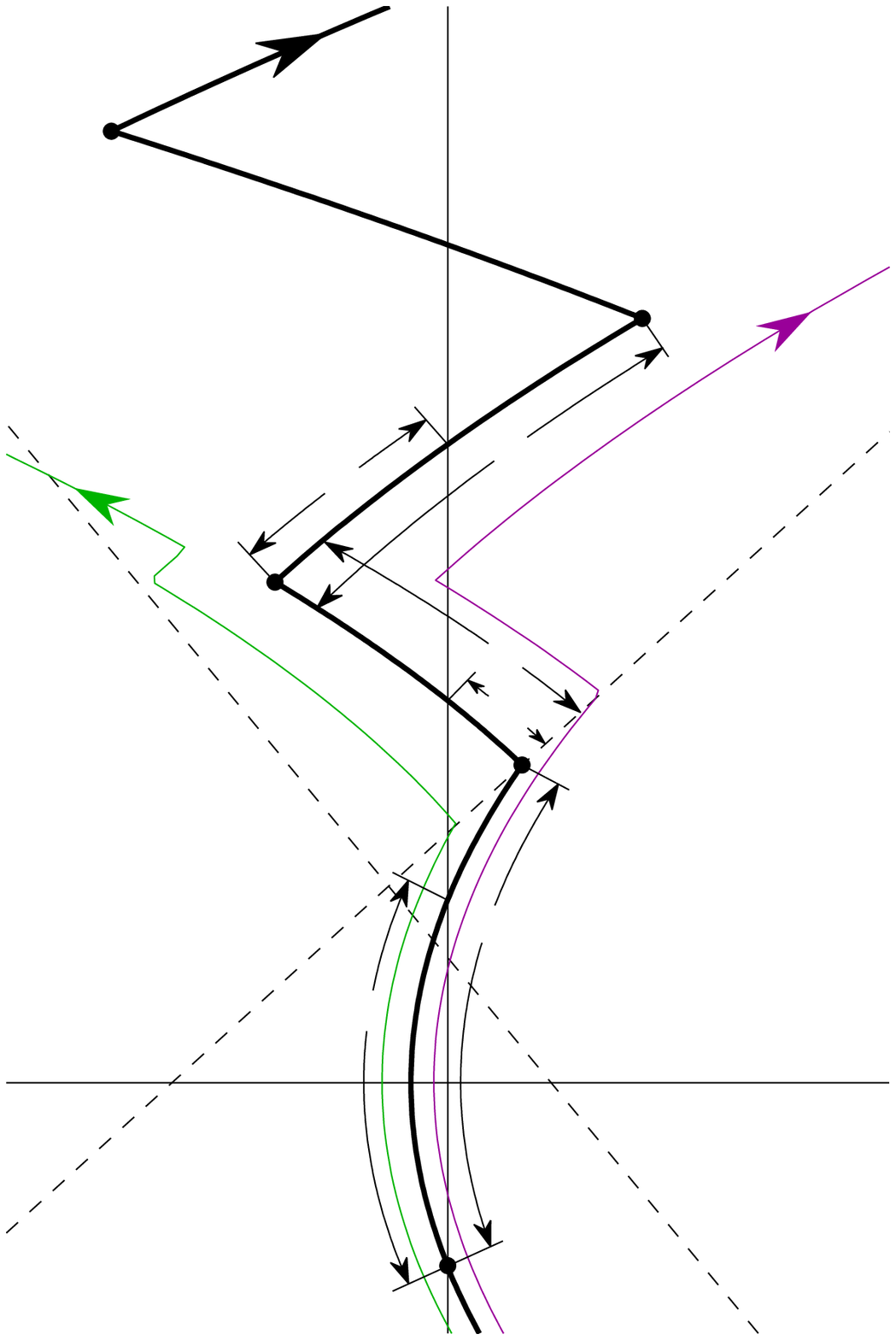}}
\put(7.6,2.35){$x$}
\put(4.06,11.6){$y$}
\put(.7,1.1){$\Gamma^-$}
\put(6.3,.7){$\Gamma^+$}
\put(4.3,.5){\scriptsize $(0,z_0)$}
\put(4.94,5.1){\scriptsize $(x_1,z_1)$}
\put(1.36,6.72){\scriptsize $(x_2,z_2)$}
\put(5.8,9.3){\scriptsize $(x_3,z_3)$}
\put(0,11.04){\scriptsize $(x_4,z_4)$}
\put(4.28,3.63){\scriptsize $\ee$}
\put(3.16,2.84){\scriptsize $T_1$}
\put(4.35,6.08){\scriptsize $S_1$}
\put(4.39,5.56){\scriptsize $T_2$}
\put(4.44,7.92){\scriptsize $S_2$}
\put(2.9,7.62){\scriptsize $T_3$}
\end{picture}
\caption{
A sketch of typical dynamics near the origin of the time-delayed system (\ref{eq:timedelay}).
\label{fig:schemTimeDelayUpper}
}
\end{center}
\end{figure}

As shown in Fig.~\ref{fig:schemTimeDelayLower},
upon passing close to the origin,
most forward orbits cease switching and either head right or head left.
However, there are two critical orbits that switch across the positive $y$-axis indefinitely,
or at least until global features of the system induce different dynamical behaviour.
As in the previous section,
we let $y_k$ for $k \ge 1$ denote the intersections of these orbits with the negative $y$-axis
at which the orbits are governed by the right half-system, see Fig.~\ref{fig:schemTimeDelayLower}.
With the definition (\ref{eq:qIdef}), we again have
\begin{equation}
q_{[y_k,y_{k-2}]} = \frac{y_{k-1} - y_k}{y_{k-2} - y_k} \;, {\rm ~for~} k = 3,5,7,\ldots \;,
\label{eq:qktd}
\end{equation}
and our goal is to calculate
\begin{equation}
\varrho_{\rm td} \equiv \lim_{\ee \to 0} q_{[y_{\rm min},y_{\rm max}]} \;,
\label{eq:varrhotimedelay}
\end{equation}
for some $-1 \ll y_{\rm min} < y_{\rm max} < 0$.
The remainder of this section consists of a derivation of formulas involving the switching points of the two critical orbits,
a recursion relation for $y_k$,
and a numerical evaluation of (\ref{eq:varrhotimedelay}).

The orbit of (\ref{eq:timedelay}) that is shown in Fig.~\ref{fig:schemTimeDelayUpper}
is located at $(0,z_0)$ at $t = 0$ and satisfies the following two assumptions.
First, it is assumed that for all $t \in (-\ee,0)$, the orbit lies to the right of $x=0$.
Consequently, until $t = \ee$, the orbit is governed by the right half-system.
Second, we assume that $z_0 < 0$ is sufficiently small that at $t = \ee$,
the orbit is located at a point $(x_1,z_1)$ with $x_1 > 0$.
By (\ref{eq:xtytR2}) we have
\begin{eqnarray}
x_1 &=& z_0 \ee + \frac{1}{2} \ee^2 + O(\ee^3) \;, \label{eq:x1} \\
z_1 &=& z_0 + \ee + O(\ee^2) \;. \label{eq:z1} 
\end{eqnarray}
The points, $(x_k,z_k)$, for $k \ge 1$, are used denote successive switching points of the orbit.
We let $T_1$ denote the time taken for the orbit to first return to $x = 0$, then
\begin{equation}
T_1 = -2 z_0 + O(\ee^2) \;.
\label{eq:T1}
\end{equation}
Similarly we let $T_k$, for each $k \ge 2$, denote the time taken for the orbit to reach $x = 0$ from each $(x_{k-1},z_{k-1})$.
Finally we let $S_k$, for each $k \ge 1$, denote the elapsed time between the
consecutive switching points $(x_k,z_k)$ and $(x_{k+1},z_{k+1})$.
Note that, $S_1$ -- the time taken for the orbit to travel between $(x_1,z_1)$ and $(x_2,z_2)$ --
is, by definition, also equal to the length of time that the orbit spends in the left half-plane
between $(0,z_0)$ and $(x_1,z_1)$, that is,
\begin{equation}
S_1 = T_1 \;.
\label{eq:S1}
\end{equation}
Similarly $S_2$ is equal to the length of time that the orbit spent previously in the right,
that is $S_2 = \ee - T_1 + T_2$,
and for $k \ge 3$, $S_k = S_{k-2} - T_{k-1} + T_k$.
From this, also by (\ref{eq:S1}), we obtain the simpler recursion relation:
\begin{equation}
S_k = \ee - S_{k-1} + T_k \;, {\rm ~for~} k \ge 2 \;.
\label{eq:Sk}
\end{equation}
Elementary use of the series solutions, (\ref{eq:xtytL2}) and (\ref{eq:xtytR2}), leads to 
\begin{equation}
\begin{array}{rcl}
T_k &=& \frac{1}{\cB} \left( \sqrt{z_{k-1}^2 + \frac{2 \cB}{\cA} x_{k-1}} - z_{k-1} \right) + O(\ee^2) \;, \\
x_k &=& x_{k-1} - \cA z_{k-1} S_{k-1} - \frac{\cA \cB}{2} S_{k-1}^2 + O(\ee^3) \;, \\
z_k &=& z_{k-1} + \cB S_{k-1} + O(\ee^2) \;,
\end{array}
\qquad {\rm for~} k = 2,4,6,\ldots \;,
\label{eq:Tkxkzkeven}
\end{equation}
and
\begin{equation}
\begin{array}{rcl}
T_k &=& \sqrt{z_{k-1}^2 - 2 x_{k-1}} - z_{k-1} + O(\ee^2) \;, \\
x_k &=& x_{k-1} + z_{k-1} S_{k-1} + \frac{1}{2} S_{k-1}^2 + O(\ee^3) \;, \\
z_k &=& z_{k-1} + S_{k-1} + O(\ee^2) \;,
\end{array}
\qquad {\rm for~} k = 3,5,7,\ldots \;.
\label{eq:Tkxkzkodd}
\end{equation}

For any values of $\cA$ and $\cB$, the above recursion relations may be used to numerically identify
the value, $z_0 = y_1(\cA,\cB)$, plotted in the inset of (\ref{fig:schemTimeDelayLower}),
for which, in the limit $\ee \to 0$, the orbit switches across the positive $y$-axis infinitely many times.
The analogous value, $\tilde{y}_1 < 0$, for the orbit that initially heads right, then switches indefinitely,
may obtained in a similar fashion.
Alternatively, by symmetry,
\begin{equation}
\tilde{y}_1(\cA,\cB) = y_1 \left( \frac{1}{\cA},\frac{1}{\cB} \right) \;.
\label{eq:tildey1}
\end{equation}
Note that it appears to be not possible to obtain an explicit expression for $y_1(\cA,\cB)$
from the recursion relations\removableFootnote{
We can reduce the above recursion relations down to just three recursion relations
valid for odd values of $k$:
\begin{equation}
\begin{split}
x_k &= x_{k-2} - \cA z_{k-2} S_{k-2} - \frac{\cA \cB}{2} S_{k-2}^2
+ (z_{k-2} + \cB S_{k-2}) \left( \ee - S_{k-2} + \frac{1}{\cB}
\left( \sqrt{z_{k-2}^2 + \frac{2 \cB}{\cA} x_{k-2}} - z_{k-2} \right) \right) \\
&+~\frac{1}{2} \left( \ee - S_{k-2} + \frac{1}{\cB}
\left( \sqrt{z_{k-2}^2 + \frac{2 \cB}{\cA} x_{k-2}} - z_{k-2} \right) \right)^2 + O(\ee^3) \;, \\
z_k &= z_{k-2} + \cB S_{k-2} + \ee - S_{k-2}
+ \frac{1}{\cB} \left( \sqrt{z_{k-2}^2 + \frac{2 \cB}{\cA} x_{k-2}} - z_{k-2} \right) + O(\ee^2) \;, \\
S_k &= S_{k-2} - \frac{1}{\cB} \left( \sqrt{z_{k-2}^2 + \frac{2 \cB}{\cA} x_{k-2}} - z_{k-2} \right) \\
&+~\sqrt{(z_{k-2} + \cB S_{k-2})^2 - 2 x_{k-2} + 2 \cA z_{k-2} S_{k-2} + \cA \cB S_{k-2}^2}
- (z_{k-2} + \cB S_{k-2}) + O(\ee^2) \;.
\end{split}
\label{eq:recursionUpperTimeDelay}
\end{equation}
The equations (\ref{eq:recursionUpperTimeDelay}) represent a three-dimensional map
and we can think of $y_1$ as the value of $z_0$ for which iterates of (\ref{eq:recursionUpperTimeDelay})
remain on a particular unstable manifold.
However this doesn't seem helpful because (\ref{eq:recursionUpperTimeDelay}) is too complicated to
invert algebraically.
If we substitute the ansatz,
$x_k = c_1 k \ee^2 + O(\ee^3)$,
$z_k = c_2 k \ee + O(\ee^2)$,
$S_k = c_3 \ee + O(\ee^2)$,
into (\ref{eq:recursionUpperTimeDelay}) and take $k \to \infty$, we find
\begin{equation}
\frac{c_1}{c_2} = \frac{\cA^2}{1+\cA+\cA^2} \;, \qquad
c_3 = \frac{1+\cA}{1+\cA+\cA^2} \;.
\end{equation}
Therefore,
\begin{equation}
\frac{S_k}{\ee} \to \frac{1+\cA}{1+\cA+\cA^2} \;,
{\rm ~as~} \ee \to 0 \;, {\rm ~for~odd~} k {\rm ~as~} k \to \infty \;.
\end{equation}
For any odd value of $k$,
we can approximate $y_1(\cA,\cB)$ by the value of $z_0$ for which, in the limit $\ee \to 0$,
$\frac{S_k}{\ee} = \frac{1+\cA}{1+\cA+\cA^2}$.
As $k \to \infty$ this approximation approaches the true value of $y_1(\cA,\cB)$.
Substituting $k = 1$ gives the crudest approximation:
$y_1 \approx -\frac{1+\cA}{2(1+\cA+\cA^2)} \ee$,
which is relatively accurate if $\cA$ is small.
}.

For backward evolution on the negative $y$-axis, illustrated in Fig.~\ref{fig:schemTimeDelayLower},
the following recursion relation may be derived from (\ref{eq:xtytL2}) and (\ref{eq:xtytR2})
(in the same manner as in the previous section)\removableFootnote{
The negative-valued times that the orbit takes to go from $y_k$ to $\Gamma^+$,
from $\Gamma^+$ to $\Gamma^-$,
and then from $\Gamma^-$ to $y_{k+2}$, are, respectively,
\begin{eqnarray}
t_1 &=& -y_k - \cA \ee - \sqrt{y_k^2 + \cA (\cA+\cB) \ee^2} + O(\ee^2) \;, \\
t_2 &=& -\frac{y_k + t_1}{\cB} - \frac{\ee}{\cA}
- \frac{1}{\cA \cB} \sqrt{ \cA^2 (y_k + t_1 - \ee \cB)^2 +
(\cA+\cB) \cB \ee^2} + O(\ee^2)\;, \\
t_3 &=& -\ee \;,
\end{eqnarray}
}\removableFootnote{
If $\cB = 0$, (\ref{eq:ykplus2td}) reduces to
$y_{k+2} = -\left( 1 + \cA + \sqrt{y_k^2/\ee^2 + \cA^2} \right) \ee + O(\ee^2)$.
}:
\begin{equation}
y_k = -\frac{1}{\cA} \left(
\cA + \cB + \sqrt{\cA^2 \left( \cA + \cB
+ \sqrt{y_{k-2}^2/\ee^2 + \cA(\cA+\cB)} \right)^2
+ (\cA+\cB) \cB} \right) \ee + O(\ee^2) \;, {\rm ~for~} k \ge 3 \;.
\label{eq:ykplus2td}
\end{equation}
Also
\begin{equation}
y_2 = -\frac{1}{\cA} \left( \cA + \cB + \sqrt{ \cA^2 \,\tilde{y}_1^2/\ee^2 + (\cA+\cB) \cB} \right) \ee + O(\ee^2) \;.
\label{eq:y2td}
\end{equation}

\begin{figure}[t!]
\begin{center}
\setlength{\unitlength}{1cm}
\begin{picture}(9.73,7)
\put(.4,0){\includegraphics[height=7cm]{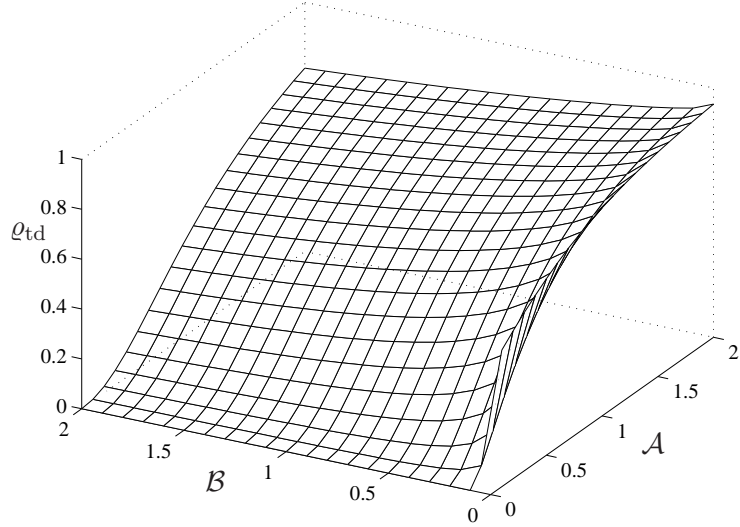}}
\put(8.4,1){\small $\cA$}
\put(2.6,.5){\small $\cB$}
\put(0,3.9){\small $\varrho_{\rm td}$}
\end{picture}
\caption{
\label{fig:ABtimedelay}
The fraction of orbits that head right after passing near the two-fold
for the time-delayed system (\ref{eq:timedelay}) in the limit $\ee \to 0$, $\varrho_{\rm td}(\cA,\cB)$,
computed numerically as described in the text.
}
\end{center}
\end{figure}

The formulas are sufficient for us to perform an efficient numerical evaluation of the fraction, $q_{[y_k,y_{k-2}]}$.
To do this, for any $\cA$ and $\cB$,
$y_1$ is computed via an iterative scheme that identifies a value of $z_0$
for which recursion according to (\ref{eq:Sk}), (\ref{eq:Tkxkzkeven}) and (\ref{eq:Tkxkzkodd}),
and starting from (\ref{eq:x1})-(\ref{eq:S1}),
gives values $(x_k,z_k)$ that lie on successively opposite sides of $x=0$
for as many values of $k$ as one chooses to compute.
Then (\ref{eq:tildey1}) is used to repeat this process and compute $\tilde{y}_1$.
Successive $y_k$ are computed via (\ref{eq:ykplus2td}) and (\ref{eq:y2td}),
and finally $q_{[y_k,y_{k-2}]}$ is evaluated by (\ref{eq:qktd}).

It is expected that
$\varrho_{\rm td} = \lim_{k \to \infty} \lim_{\ee \to 0} q_{[y_k,y_{k-2}]}$,
since this identity holds for the hysteretic perturbation of the previous section
(it is an elementary consequence of Lemma \ref{le:qI}).
However, since here an exact expression for $y_k$ is not known,
and in view of the complexity in proving Lemma \ref{le:qI},
a verification of this claim is left for future studies.
Fig.~\ref{fig:ABtimedelay} shows a plot of $\varrho_{\rm td}$ approximated by evaluating $q_{[y_k,y_{k-2}]}$
at suitably large $k$ and small $\ee$ via the iterative scheme described above\removableFootnote{
A little algebra reveals that $\lim_{\cA \to 0} \varrho_{\rm td}(\cA,\cB) = 0$,
as indicated in Fig.~\ref{fig:ABtimedelay}.
First, as $\cA \to 0$, it follows from (\ref{eq:ykplus2td}) and (\ref{eq:y2td})
that $y_2 \approx y_3 \approx - \left( 1 + \frac{2 \cB}{\cA} \right) \ee$,
which are therefore large relative to $y_1$.
Consequently $p_1$ is small and, in the limit $\cA \to 0$, $p_1 \to 0$ .
Therefore, as $\cA \to 0$, each $q_k \to 0$ and thus $\varrho_{\rm td} \to 0$.
}.

\section{Noise}
\label{sec:NOISE}
\setcounter{equation}{0}

This section concerns the following formulation of (\ref{eq:ode}) in the presence of noise:
\begin{equation}
\left[ \begin{array}{c} dx(t) \\ dy(t) \end{array} \right] =
\left\{ \begin{array}{lc}
\left[ \begin{array}{c} -\cA y + O(|x|) + O(2) \\ \cB + O(1) \end{array} \right] \;, & x < 0 \\
\left[ \begin{array}{c} y + O(|x|) + O(2) \\ 1 + O(1) \end{array} \right] \;, & x > 0
\end{array} \right\} \,dt +
\ee \left[ \begin{array}{cc}
d_{11} & d_{12} \\
d_{21} & d_{22}
\end{array} \right]
\left[ \begin{array}{c} dW_1(t) \\ dW_2(t) \end{array} \right] \;.
\label{eq:sde}
\end{equation}
Here $\ee$ denotes the overall noise amplitude
and $d_{ij}$ are constants that control the relative strength of the two-dimensional noise
in the $x$ and $y$ directions.
Alternative formulations include coloured noise, spatially dependent noise, and multiplicative or parametric noise,
each of which may be more suitable in certain applications.

\subsection{General aspects of the zero-noise limit}

Equation (\ref{eq:sde}) is an example of a multi-dimensional stochastic differential equation,
\begin{equation}
d\bx(t) = \phi(\bx) \,dt + \ee D \,d\bW(t) \;,
\label{eq:sdeGen}
\end{equation}
for which the drift, $\phi$, may be discontinuous.
If $\phi$ is measurable and bounded\removableFootnote{
I think there still exists a strong solution if the vector field grows at most linearly
(i.e.~there exists $M > 0$ such that $||\phi(\bx)|| \le M (1 + ||\bx||)$)
which is sufficient to ensure that solutions to the deterministic solution do not
become infinite in finite time.
}
and the diffusion matrix, $D$, is non-singular\removableFootnote{
Really we need $D D^{\sf T}$ to be positive definite (that is, $v^{\sf T} D D^{\sf T} v > 0$, for all $v \ne 0$).
But recall $D D^{\sf T}$ is positive definite if and only if $\det(D) \ne 0$
(simply because $v^{\sf T} D D^{\sf T} v = 0$ if and only if $v$ is an eigenvector of $D^{\sf T}$ with eigenvalue $0$).
},
then (\ref{eq:sdeGen}) has a unique strong solution\removableFootnote{
Recall, a {\em strong solution}
is a function $\bx(t) = \mathcal{F}(t,\bW(s),0<s<t)$
(i.e.~given explicitly in terms of $\bW$)
whereas a {\em weak solution} is a process with the same distribution as $\bx(t)$,
but not given in terms of $\bW$ \cite{Ok03}.
The Tanaka equation, $dx(t) = {\rm sgn}(x(t)) \,dW(t)$,
has a weak solution but no strong solution \cite{Ok03}.
}
from any initial point \cite{Fl11,KrRo05,PrSh98,StVa69}\removableFootnote{
See also \cite{Ve81,Co71,Co72,LeRi85,MoNi12}.
}.
In following discussion, $p_\ee(\bx,t;\bx_0)$ denotes the probability density function of $\bx(t)$, given $\bx(0) = \bx_0$.

With $\ee = 0$, (\ref{eq:sdeGen}) reduces to the ODE
\begin{equation}
\dot{\bx} = \phi(\bx) \;.
\label{eq:odeGen}
\end{equation}
If $\phi$ is smooth, (\ref{eq:odeGen}) has a unique solution, call it $\bx_{\rm det}(t;\bx_0)$.
Straight-forward asymptotic arguments are sufficient to demonstrate weak convergence:
$p_\ee(\bx,t;\bx_0) \to \delta(\bx - \bx_{\rm det}(t;\bx_0))$ as $\ee \to 0$,
where $\delta$ denotes the Dirac-delta function \cite{FrWe12,Ga09}.
If $\phi$ is continuous but non-differentiable, (\ref{eq:odeGen}) may have multiple solutions
and $p_\ee$ may converge to a weighted sum of delta functions centred at some of these solutions.
In particular we may have
\begin{equation}
p_\ee(\bx,t;\bx_0) \to (1-\varrho) \delta(\bx - \bx_1(t;\bx_0)) + \varrho \delta(\bx - \bx_2(t;\bx_0)) \;,
\label{eq:sumtwodeltafns}
\end{equation}
as $\ee \to 0$,
where $\bx_1(t;\bx_0)$ and $\bx_2(t;\bx_0)$ are two different solutions to (\ref{eq:odeGen}).

In the short article \cite{Ve83}, stochastic calculus methods are used to
prove (\ref{eq:sumtwodeltafns}) in one dimension 
when $\phi$ is continuous, odd, nondecreasing, convex on $[0,\infty)$,
and with a finite value of $\lim_{\Delta \to 0} \int_0^\Delta \frac{1}{\phi(x)} \,dx$, when $x_0 = 0$.
Here, $\varrho = \frac{1}{2}$ by symmetry, and of the uncountably many solutions to (\ref{eq:odeGen}),
the $x_1$ and $x_2$ that appear in (\ref{eq:sumtwodeltafns})
are the two solutions for which $x(t) \ne 0$ for all $t > 0$.
These are called the {\em minimal} and {\em maximal} solutions
because, for any other solution, $x(t)$, for all $t > 0$, $x_1(t) < x(t) < x_2(t)$.
In \cite{BaBa82}, similar results are obtained when $\phi$ is not odd
and the value of $\varrho$ is derived using first passage concepts.
Similar results for multiple dimensions have recently been described in \cite{Zh12}.

For discontinuous $\phi$, it is necessary to consider generalised notions of a solution to (\ref{eq:odeGen}) \cite{Co08c}.
In \cite{BuOu09} was is shown that solutions to (\ref{eq:sdeGen})
may limit only to Filippov solutions\removableFootnote{
A {\em Filippov solution} to (\ref{eq:odeGen}) is an {\em absolutely continuous}
(absolute continuity is a stronger condition than uniform continuity;
an example of a uniformly continuous function that is not absolutely continuous is the Cantor function;
absolutely continuous functions are differentiable almost everywhere)
function, $\bx(t)$, that satisfies the {\em differential inclusion} \cite{De92,Sm00},
\begin{equation}
\dot{\bx} \in \bigcap_{\mu(\Sigma) = 0} \bigcap_{\Delta > 0}
\,{\rm co} \left[ \phi \left( (\bx + \Delta B) \setminus \Sigma \right) \right] \;,
\label{eq:differentialInclusion}
\end{equation}
almost everywhere.
Let us explain (\ref{eq:differentialInclusion}).
$B = \left\{ \bx ~\big| ||\bx|| < 1 \right\}$ denotes the unit ball in $\mathbb{R}^N$,
and $\mu$ denotes Lebesgue measure.
Thus $(\bx + \Delta B) \setminus \Sigma$ is the set of all points within a distance $\Delta$ of $\bx$,
except those on some measure zero set $\Sigma$.
Also, ${\rm co}(\Omega) \equiv \left\{ (1-q) \bx_1 + q \bx_2 ~\big|~
0 \le q \le 1 \;, \bx_1, \bx_2 \in \Omega \right\}$, denotes the closed convex hull of a set $\Omega$.
By taking the intersection over all $\Delta > 0$ we are essentially taking the limit $\Delta \to 0$,
i.e.~(\ref{eq:differentialInclusion}) is a local definition.
By neglecting measure zero sets, $\Sigma$, and taking the intersection over all such $\Sigma$,
(\ref{eq:differentialInclusion}) is not influenced by, for instance, values of $\phi$ at a switching manifold.
}
of (\ref{eq:odeGen}) as the noise amplitude is taken to zero\removableFootnote{
Similar results have been achieved using Malliavin calculus \cite{MeMe11},
and related limit theorems are described in \cite{Ma98}, but I don't understand these results.
}.
However, in the case that there are multiple solutions,
this does not tell us which Filippov solutions contribute to the limiting density.
On stable sliding regions there is a unique Filippov solution,
and if sliding occurs for the duration of time interval, $[0,t]$,
it may be shown\removableFootnote{
Assuming the asymptotic expression is valid (which we haven't proved),
it should be straight-forward to show that higher order terms may be neglected
and the difference between the leading order term and the delta-function goes to zero in some norm.
}
directly that $p_\ee(\bx,t;\bx_0)$ limits to a delta-function centred at the Filippov solution as $\ee \to 0$,
via an asymptotic expansion of $p$ \cite{SiKu13b}.

If $\phi$ is discontinuous, then, to the author's knowledge,
only when $\phi$ is one-dimensional and piecewise-constant is an exact expression for $p_\ee(\bx,t;\bx_0)$ known.
Specifically, for any constants $\phi^\pm \in \mathbb{R}$,
the transitional probability density function of
\begin{equation}
dx(t) = \left\{ \begin{array}{lc}
\phi^- \;, & x < 0 \\
\phi^+ \;, & x > 0
\end{array} \right\} \,dt + \ee \,dW(t) \;, \qquad x(0) = 0 \;,
\label{eq:sdePWC}
\end{equation}
is given by
\begin{equation}
p_\ee(x,t|0) = \left\{ \begin{array}{lc}
\frac{2}{\ee^2}
\,{\rm e}^{\frac{2 \phi^- x}{\ee^2}}
\int_0^\infty \int_0^t
\frac{b-x}{\sqrt{2 \pi \ee^2 (t-u)^3}}
\,{\rm e}^{\frac{-(b - x + \phi^- (t-u))^2}{2 \ee^2 (t-u)}}
\frac{b}{\sqrt{2 \pi \ee^2 u^3}}
\,{\rm e}^{\frac{-(b - \phi^+ u)^2}{2 \ee^2 u}}
\,du \,db \;, & x \le 0 \\
\frac{2}{\ee^2}
\,{\rm e}^{\frac{2 \phi^+ x}{\ee^2}}
\int_0^\infty \int_0^t
\frac{b+x}{\sqrt{2 \pi \ee^2 (t-u)^3}}
\,{\rm e}^{\frac{-(b + x + \phi^+ (t-u))^2}{2 \ee^2 (t-u)}}
\frac{b}{\sqrt{2 \pi \ee^2 u^3}}
\,{\rm e}^{\frac{-(b - \phi^- u)^2}{2 \ee^2 u}}
\,du \,db \;, & x \ge 0
\end{array} \right. \;,
\label{eq:KaSh84}
\end{equation}
see \cite{KaSh84,KaSh91}.
Moreover, if $\phi^- < 0 < \phi^+$, then multiple Filippov solutions emanate from $x=0$,
and by using Laplace's method \cite{BeOr99} to evaluate (\ref{eq:KaSh84}) asymptotically\removableFootnote{
We do this here for $x > 0$ with $a_L = \phi^-$ and $a_R = -\phi^+$.
The result for $x < 0$ follows by symmetry.

As shown in \cite{KaSh84,KaSh91}, for $x > 0$,
\begin{equation}
p_\ee(x,t|0) = \frac{2}{\ee^2}
\,{\rm e}^{\frac{-2 a_R x}{\ee^2}}
\int_0^\infty \int_0^t
\frac{b+x}{\sqrt{2 \pi \ee^2 (t-u)^3}}
\,{\rm e}^{\frac{-(b+x-a_R(t-u))^2}{2 \ee^2 (t-u)}}
\frac{b}{\sqrt{2 \pi \ee^2 u^3}}
\,{\rm e}^{\frac{-(b-a_L u)^2}{2 \ee^2 u}}
\,du \,db \;.
\end{equation}
We may write this as
\begin{equation}
p_\ee(x,t|0) = \frac{1}{\pi \ee^4} \int_0^\infty \int_0^t
g(b,u) \,{\rm e}^{\frac{-1}{\ee^2} f(b,u)} \,du \,db \;,
\end{equation}
where
\begin{eqnarray}
f(b,u) &=& 2 a_R x + \frac{(b+x-a_R(t-u))^2}{2 (t-u)} + \frac{(b-a_L u)^2}{2 u} \;, \\
g(b,u) &=& \frac{b(b+x)}{u^{\frac{3}{2}} (t-u)^{\frac{3}{2}}} \;.
\end{eqnarray}
We apply Laplace's method for the integral in $u$
about the value, $u^*(b) = k b + O(b^2)$, at which $f(b,u)$ is minimal (defined by $f_u(b,u^*(b)) \equiv 0$).
We have
\begin{eqnarray}
f(b,u^*(b)) &=& \frac{(x+a_R t)^2}{2 t} + l b + O(b^2) \;, \\
f_{uu}(b,u^*(b)) &=& \frac{1}{b} \left( \frac{1}{k^3} + O(b) \right) \;, \\
g(b,u^*(b)) &=& \frac{1}{\sqrt{b}} \left( \frac{x}{t^{\frac{3}{2}} k^{\frac{3}{2}}} + O(b) \right) \;,
\end{eqnarray}
where
\begin{equation}
k = \frac{1}{\sqrt{|a_L^2 - a_R^2 + \frac{x^2}{t^2}|}} \;, \qquad
l = \frac{(x-a_R t)^2 k}{t^2} + \frac{(x-a_R t)(1+a_R k)}{t} + \frac{(1-a_L k)^2}{2 k} \;.
\end{equation}
Upon substituting, $w = \frac{u - u^*(b)}{\ee}$, integrating and simplifying we arrive at
\begin{equation}
p_\ee(x,t|0) = \left( \frac{\sqrt{2} x}{\sqrt{\pi} \ee^3 t^{\frac{3}{2}}}
\int_0^\infty \left( 1 + O(b) \right)
\,{\rm e}^{\frac{-1}{\ee^2}
\left( \frac{(x+a_R t)^2}{2 t} + l b + O(b^2) \right)} \,db \right)
\left( 1 + O(\ee) \right) \;.
\end{equation}
A second application of Laplace's method,
this time to the integral in $b$ about zero (using, say $q = \frac{b}{\ee^2}$), yields
\begin{equation}
p_\ee(x,t|0) = \frac{\sqrt{2} x}{l \sqrt{\pi} \ee t^{\frac{3}{2}}}
\,{\rm e}^{\frac{-(x+a_R t)^2}{2 \ee^2 t}} + O(\ee^0) \;.
\end{equation}
Finally, expanding $x$ about $-a_R t$ gives,
\begin{equation}
p_\ee(x,t|0) = \frac{a_R}{a_L+a_R} \frac{1}{2 \pi \ee^2 t}
\,{\rm e}^{\frac{-(x+a_R t)^2}{2 \ee^2 t}} + O(|x + a_R t|) + O(\ee^0) \;,
\end{equation}
where, in particular, $l = -2(a_L+a_R) + O(|x+a_R t|)$.
We conclude that,
\begin{equation}
p_\ee(x,t|0) \to \frac{a_R}{a_L+a_R} \delta(x + a_R t) \;,
\end{equation}
in the $L^1$-norm, as $\ee \to 0$.
},
it may be shown that
\begin{equation}
p_\ee(x,t;0) \to \frac{-\phi^-}{\phi^+ - \phi^-} \,\delta(x - \phi^- t) + \frac{\phi^+}{\phi^+ - \phi^-} \,\delta(x - \phi^+ t) \;.
\end{equation}
This matches the results for continuous $\phi$, in that $\phi^- t$ and $\phi^+ t$
are the minimal and maximal Filippov solutions to
(\ref{eq:sdePWC}) with $\ee = 0$.
Large deviations for the case $(\phi^-,\phi^+) = (-1,1)$ and generalisations are described in \cite{GrHe01}.

\begin{figure}[b!]
\begin{center}
\setlength{\unitlength}{1cm}
\begin{picture}(5.6,7)
\put(0,0){\includegraphics[height=7cm]{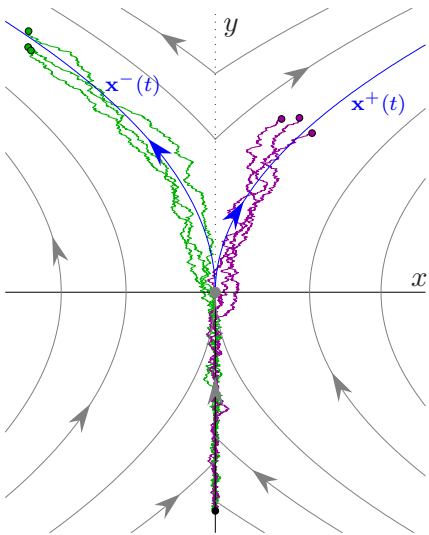}}
\put(5.4,3.3){\small $x$}
\put(2.9,6.7){\small $y$}
\put(1.33,5.9){\scriptsize \color{blue} $\bx^-(t)$}
\put(4.6,5.65){\scriptsize \color{blue} $\bx^+(t)$}
\end{picture}
\caption{
Six typical sample paths of (\ref{eq:sde}).
\label{fig:schemNoise}
}
\end{center}
\end{figure}

\subsection{Asymptotics of a first passage representation}

Our interest is in forward evolution of (\ref{eq:sde}) from points on the negative $y$-axis.
For small $\ee > 0$, upon passing near the two-fold,
sample paths follow close to either $\bx^+(t)$ or $\bx^-(t)$ with high probability.
This is illustrated in Fig.~\ref{fig:schemNoise}.
Here the probability that a sample path will eventually head right or head left
is determined via first passage methods for arbitrarily small $\ee > 0$.

We restrict our attention to an $\ee$-independent rectangle,
$[-x^*,x^*] \times [-y^*,y^*]$, in the $(x,y)$-plane,
and assume $x^*,y^* > 0$ are sufficiently small that global features of the dynamics may be ignored.
For any initial point, $(x(0),y(0)) = (x_0,y_0)$, inside the rectangle,
we let $Q_\ee(x_0,y_0)$ denote the probability that first passage by (\ref{eq:sde})
to the boundary of the rectangle occurs at a point with $x > 0$.
We then define
\begin{equation}
\varrho_{\rm no} \equiv \lim_{\ee \to 0} Q_\ee(0,y_0;x^*,y^*) \;,
\label{eq:varrhonoise0}
\end{equation}
for, $-y^* < y_0 < 0$,
as this represents the probability that forward evolution from the negative $y$-axis heads right
after passing by the origin in the zero-noise limit.
Below it will become apparent that this limiting probability is independent of the values of $x^*$, $y^*$ and $y_0$
as implied by (\ref{eq:varrhonoise0}).

As described in stochastic differential equation texts, such as \cite{Sc10}, $Q_\ee$ is the $C^1$
solution to the elliptic boundary value problem comprised of the backward Fokker-Planck equation for (\ref{eq:sde}):
\begin{eqnarray}
0 &=& \left\{ \begin{array}{lc}
\left( -\cA y_0 + O(|x_0|) + O(2) \right) \frac{\partial Q_\ee}{\partial x_0} +
\left( \cB + O(1) \right) \frac{\partial Q_\ee}{\partial y_0} \;, & x_0 < 0 \\
\left( y_0 + O(|x_0|) + O(2) \right) \frac{\partial Q_\ee}{\partial x_0} +
\left( 1 + O(1) \right) \frac{\partial Q_\ee}{\partial y_0} \;, & x_0 > 0
\end{array} \right\} \nonumber \\
&&~+
\frac{\ee}{2} (d_{11}^2 + d_{12}^2) \frac{\partial^2 Q_\ee}{\partial x_0^2} +
\ee (d_{11} d_{21} + d_{12} d_{22}) \frac{\partial^2 Q_\ee}{\partial x_0 \partial y_0} +
\frac{\ee}{2} (d_{21}^2 + d_{22}^2) \frac{\partial^2 Q_\ee}{\partial y_0^2} \;,
\label{eq:Qeepde}
\end{eqnarray}
and the boundary conditions
\begin{equation}
Q_\ee(x_0,\pm y^*) = \chi_{[0,\infty)}(x_0) \;, \qquad
Q_\ee(-x^*,y_0) = 0 \;, \qquad
Q_\ee(x^*,y_0) = 1 \;.
\label{eq:bcs0}
\end{equation}
The boundary conditions follow from the observation that first passage occurs in zero time
for any point on the boundary of the rectangle.
Therefore, as stated in (\ref{eq:bcs0}), on the boundary, $Q_\ee = 1$ if $x_0 > 0$, and $Q_\ee = 0$ otherwise.

A description of $Q_\ee$ for small $\ee$ may be obtained from an asymptotic expansion of (\ref{eq:Qeepde})-(\ref{eq:bcs0}).
Since only $\lim_{\ee \to 0} Q_\ee$ is of interest,
only the leading order term in this expansion is required here.
A formal justification of the expansion 
is a difficult analytical task deferred for future work\removableFootnote{
A proof will probably require the maximum principle \cite{BeOr99,VaBu95}.
I can do such a proof for the analogous problem in one dimension,
but am stuck searching for an analogous proof for unstable sliding in two dimensions
(which should be a slightly simpler problem than the one here).
Note, often asymptotic solutions of this nature are correct in the limit $\ee \to 0$,
but comprised of a divergent series \cite{Ho95}.
}.

Here the scaling laws of the asymptotic expansion are obtained in an intuitive fashion
via a blow-up of the dynamics of (\ref{eq:sde}) about the origin.
Specifically, with the general scaling,
\begin{equation}
\tilde{x} = \frac{x}{\ee^{\lambda_1}} \;, \qquad
\tilde{y} = \frac{y}{\ee^{\lambda_2}} \;, \qquad
\tilde{t} = \frac{t}{\ee^{\lambda_3}} \;,
\label{eq:scaling}
\end{equation}
where each $\lambda_i > 0$, (\ref{eq:sde}) becomes
\begin{eqnarray}
\left[ \begin{array}{c} d\tilde{x}(\tilde{t}) \\ d\tilde{y}(\tilde{t}) \end{array} \right] &=&
\left\{ \begin{array}{lc}
\left[ \begin{array}{c}
-\ee^{-\lambda_1 + \lambda_2 + \lambda_3} \cA \tilde{y}(\tilde{t})
+ O \left( \ee^{-\lambda_1 + \lambda_3 + {\rm min}(\lambda_1,2 \lambda_2)} \right) \\
\ee^{-\lambda_2 + \lambda_3} \cB
+ O \left( \ee^{-\lambda_2 + \lambda_3 + {\rm min}(\lambda_1,\lambda_2)} \right)
\end{array} \right] \;, & \tilde{x} < 0 \\
\left[ \begin{array}{c} 
\ee^{-\lambda_1 + \lambda_2 + \lambda_3} \tilde{y}(\tilde{t})
+ O \left( \ee^{-\lambda_1 + \lambda_3 + {\rm min}(\lambda_1,2 \lambda_2)} \right) \\
\ee^{-\lambda_2 + \lambda_3}
+ O \left( \ee^{-\lambda_2 + \lambda_3 + {\rm min}(\lambda_1,\lambda_2)} \right)
\end{array} \right] \;, & \tilde{x} > 0
\end{array} \right\} \,d\tilde{t} \nonumber \\
&&+~\left[ \begin{array}{cc}
\ee^{1 - \lambda_1 + \frac{\lambda_3}{2}} d_{11} & \ee^{1 - \lambda_1 + \frac{\lambda_3}{2}} d_{12} \\
\ee^{1 - \lambda_2 + \frac{\lambda_3}{2}} d_{21} & \ee^{1 - \lambda_2 + \frac{\lambda_3}{2}} d_{22}
\end{array} \right]
\left[ \begin{array}{c} dW_1(\tilde{t}) \\ dW_2(\tilde{t}) \end{array} \right] \;.
\label{eq:sde2}
\end{eqnarray}
The appropriate values of $\lambda_i$ are those for which three terms in (\ref{eq:sde2}) are $O(1)$
and all other terms are of higher order.
This gives two possibilities\removableFootnote{
Assuming $\lambda_i > 0$ for each $i$,
the exponents in $\ee$ of the four potentially leading order terms are
\begin{equation}
-\lambda_1 + \lambda_2 + \lambda_3 \;, \qquad
-\lambda_2 + \lambda_3 \;, \qquad
1 - \lambda_1 + \frac{\lambda_3}{2} \;, \qquad
1 - \lambda_2 + \frac{\lambda_3}{2} \;.
\end{equation}
The first three quantities are zero when
$(\lambda_1,\lambda_2,\lambda_3) = \left( \frac{4}{3}, \frac{2}{3}, \frac{2}{3} \right)$
and then the fourth quantity is $\frac{2}{3}$.\\
The first two and last quantities are zero when
$(\lambda_1,\lambda_2,\lambda_3) = (4,2,2)$
and then the third quantity is $-2$.\\
The first and last two quantities are zero when
$(\lambda_1,\lambda_2,\lambda_3) = (1,1,0)$
and then the second quantity is $-1$.\\
The last three quantities are zero when
$(\lambda_1,\lambda_2,\lambda_3) = (2,2,2)$
and then the first quantity is $2$.
}.
First we may have $\lambda_1 = \lambda_2 = \lambda_3 = 2$,
but this is merely the scaling that in general sets drift and diffusion to the same order
and thus covers a time-scale too short to describe stochastic departure
from the discontinuity surface. 
Hence, we use the second possibility, which is
\begin{equation}
\lambda_1 = \frac{4}{3} \;, \qquad
\lambda_2 = \frac{2}{3} \;, \qquad
\lambda_3 = \frac{2}{3} \;.
\label{eq:lambdas}
\end{equation}
Since (\ref{eq:lambdas}) scales $x$ to a higher degree than $y$, noise in $y$ is of higher order than noise in $x$.
By absorbing the magnitude of the noise in $x$ into $\ee$, we may assume
\begin{equation}
d_{11}^2 + d_{12}^2 = 1 \;,
\end{equation}
so that the noise in $x$ is normalised.
In the limit, $\ee \to 0$, we define
\begin{equation}
\tilde{x} \to u \;, \qquad
\tilde{y} \to v \;, \qquad
\tilde{t} \to s \;,
\label{eq:uvs}
\end{equation}
with which the stochastic differential equation, (\ref{eq:sde2}), takes the reduced form
\begin{equation}
\left[ \begin{array}{c} du(s) \\ dv(s) \end{array} \right] =
\left\{ \begin{array}{lc}
\left[ \begin{array}{c} -\cA v \\ \cB \end{array} \right] \;, & u < 0 \\
\left[ \begin{array}{c} v \\ 1 \end{array} \right] \;, & u > 0
\end{array} \right\} \,ds +
\left[ \begin{array}{c} dW(s) \\ 0 \end{array} \right] \;.
\label{eq:sdeReduced}
\end{equation}

We now return to the boundary value problem for the probability, $Q_\ee$.
In terms of the limiting scaled variables, (\ref{eq:uvs}),
the leading order term for the asymptotic expansion of $Q_\ee$, call it $Q$,
satisfies the backward Fokker-Planck equation of (\ref{eq:sdeReduced}):
\begin{equation}
0 = \left\{ \begin{array}{lc}
-\cA v_0 \frac{\partial Q}{\partial u_0} +
\cB \frac{\partial Q}{\partial v_0} +
\frac{1}{2} \frac{\partial^2 Q}{\partial u_0^2} \;, & u_0 < 0 \\
v_0 \frac{\partial Q}{\partial u_0} +
\frac{\partial Q}{\partial v_0} +
\frac{1}{2} \frac{\partial^2 Q}{\partial u_0^2} \;, & u_0 > 0
\end{array} \right. \;.
\label{eq:Qpde}
\end{equation}
Since the boundary conditions for $Q_\ee$ (\ref{eq:bcs0}) are at the boundary of a fixed rectangle,
the boundary conditions for $Q$ are at infinity.
However, by the nature of the reduction, (\ref{eq:Qpde}) is parabolic and
consequently we cannot impose a boundary condition on $Q$ in the limit $v_0 \to -\infty$.
This is not problematic because
in the $v$-direction (\ref{eq:sdeReduced}) has positive drift and no noise,
thus for sample paths to (\ref{eq:sdeReduced}),
$v(t) \to \infty$ as $t \to \infty$ with probability $1$.
$Q$ satisfies the remaining three boundary conditions,
\begin{equation}
\lim_{v_0 \to \infty} Q(u_0,v_0) = \chi_{[0,\infty)}(u_0) \;, \qquad
\lim_{u_0 \to -\infty} Q(u_0,v_0) = 0 \;, \qquad
\lim_{u_0 \to \infty} Q(u_0,v_0) = 1 \;.
\label{eq:bcs}
\end{equation}
The goal is to determine $\varrho_{\rm no}$ (\ref{eq:varrhonoise0}),
which may now be written as $\varrho_{\rm no} = \lim_{v_0 \to -\infty} Q(0,v_0)$.

\subsection{The limiting probability, $\varrho_{\rm no}$}

\begin{figure}[t!]
\begin{center}
\setlength{\unitlength}{1cm}
\begin{picture}(9.33,7.4)
\put(0,0){\includegraphics[height=7cm]{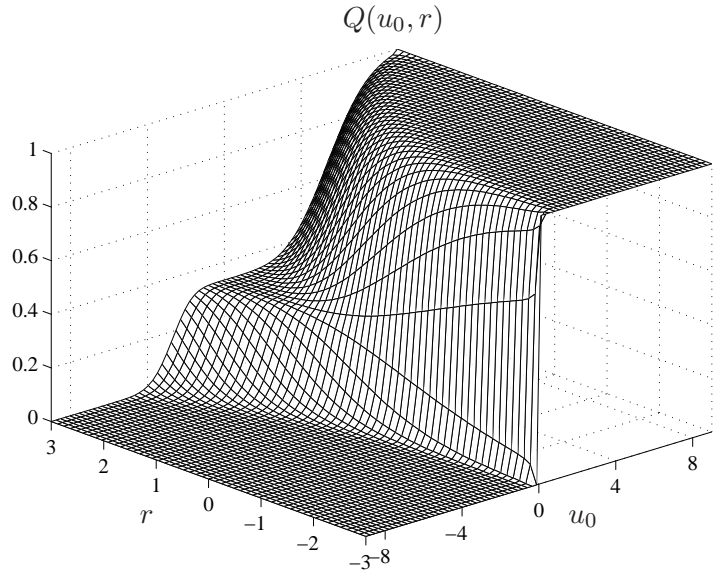}}
\put(7.4,.7){\small $u_0$}
\put(1.7,.7){\small $r$}
\put(4.4,7.3){\small $Q(u_0,r)$}
\end{picture}
\caption{
The function $Q(u_0,r)$ when $\cA = 2$ and $\cB = 4$ 
as computed by numerically solving the boundary value problem (\ref{eq:Qpde2})-(\ref{eq:bcs2}).
For these values of $\cA$ and $\cB$,
by (\ref{eq:varrhonoise}) and (\ref{eq:sV1}) we have
$Q(0,r) \to \frac{1}{3}$ as $r \to \infty$.
\label{fig:drawQ}
}
\end{center}
\end{figure}

\begin{figure}[t!]
\begin{center}
\setlength{\unitlength}{1cm}
\begin{picture}(9.73,7)
\put(.4,0){\includegraphics[height=7cm]{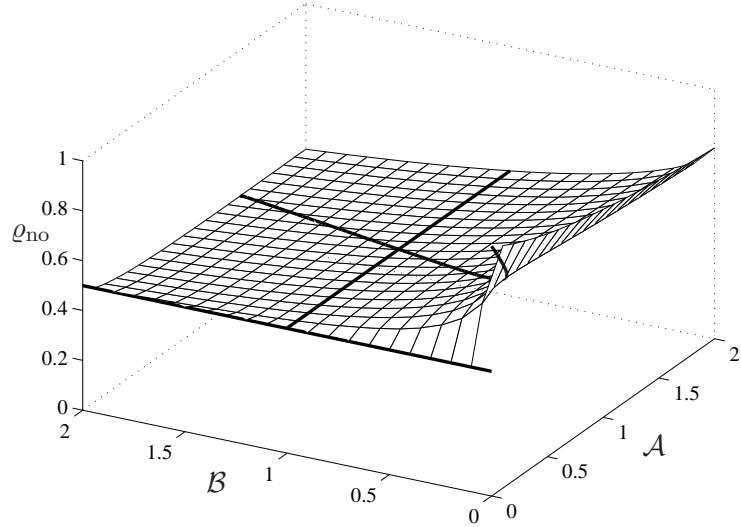}}
\put(8.4,1){\small $\cA$}
\put(2.6,.5){\small $\cB$}
\put(0,3.9){\small $\varrho_{\rm no}$}
\end{picture}
\caption{
A plot of $\varrho_{\rm no}(A,B)$ computed by numerically solving
the boundary value problem (\ref{eq:Qpde2})-(\ref{eq:bcs2}) and
evaluating $Q(0,r)$ at a large value of $r$ to approximate (\ref{eq:varrhonoise}).
The solid curves are the theoretical values of $\varrho_{\rm no}(A,B)$ as given by Lemma \ref{le:specialValues}.
\label{fig:ABnoise}
}
\end{center}
\end{figure}

The partial differential equation, (\ref{eq:Qpde}),
takes a more familiar form if we introduce the pseudo-time
\begin{equation}
r = -v_0 \;.
\end{equation}
In terms of $r$, the parabolic boundary value problem is
\begin{equation}
\frac{\partial Q}{\partial r} = \left\{ \begin{array}{lc}
\frac{\cA}{\cB} r \frac{\partial Q}{\partial u_0} +
\frac{1}{2 \cB} \frac{\partial^2 Q}{\partial u_0^2} \;, & u_0 < 0 \\
-r \frac{\partial Q}{\partial u_0} +
\frac{1}{2} \frac{\partial^2 Q}{\partial u_0^2} \;, & u_0 > 0
\end{array} \right. \;,
\label{eq:Qpde2}
\end{equation}
\begin{equation}
\lim_{r \to -\infty} Q(u_0,r) = \chi_{[0,\infty)}(u_0) \;, \qquad
\lim_{u_0 \to -\infty} Q(u_0,r) = 0 \;, \qquad
\lim_{u_0 \to \infty} Q(u_0,r) = 1 \;.
\label{eq:bcs2}
\end{equation}
Reinterpreting $Q$ as a function of $u_0$ and $r$, we have 
\begin{equation}
\varrho_{\rm no} = \varrho_{\rm no}(\cA,\cB) = \lim_{r \to \infty} Q(0,r) \;.
\label{eq:varrhonoise}
\end{equation}
Due to the discontinuity at $u_0 = 0$ and the presence of $r$ on the right hand side of (\ref{eq:Qpde2}),
it appears to be not possible to obtain
an explicit expression for $Q$ from (\ref{eq:Qpde2})-(\ref{eq:bcs2})\removableFootnote{
When $\cA = \cB = 1$, the stochastic differential equation reduces to 1d
for which the solution is given by half of Knessl's solution on each side of zero.
}.
Fig.~\ref{fig:drawQ} shows a typical plot of $Q(u_0,r)$ as computed by finite difference approximations.
Fig.~\ref{fig:ABnoise} shows a plot of $\varrho_{\rm no}(\cA,\cB)$ computed by then
approximating the limit (\ref{eq:varrhonoise}) at a suitably large value of $r$.
The following lemma provides the value of $\varrho_{\rm no}$
at some special values of $\cA$ and $\cB$:
\begin{lemma}
The limiting value, (\ref{eq:varrhonoise}), of the parabolic boundary value problem,
(\ref{eq:Qpde2})-(\ref{eq:bcs2}), satisfies 
\begin{eqnarray}
\varrho_{\rm no}(\cA,\cA^2) &=& \frac{1}{1+\cA} \;, \label{eq:sV1} \\
\varrho_{\rm no}(\cA,1) &=& \frac{1}{2} \;, \label{eq:sV2} \\
\lim_{\cA \to 0} \varrho_{\rm no}(\cA,\cB) &=& \frac{1}{2} \;. \label{eq:sV3}
\end{eqnarray}
\label{le:specialValues}
\end{lemma}
These values are highlighted in Fig.~\ref{fig:ABnoise}.
A derivation of (\ref{eq:sV1})-(\ref{eq:sV3}) is given in Appendices \ref{sec:SV1}-\ref{sec:SV3}.
Since $\cA$ and $\cB$ may take any positive values,
it follows from (\ref{eq:sV1}) that all values in $(0,1)$ are possible for $\varrho_{\rm no}$.

Finally, consider the transitional probability density function for (\ref{eq:sde}), call it $p_\ee(\bx,t;\bx_0)$.
Write $\bx_0 = (x_0,y_0)$, and suppose $x_0 = 0$ and $-1 \ll y_0 < 0$.
Let $t_{\rm slide}$ denote the time taken for the deterministic sliding trajectory of (\ref{eq:ode})
to travel from $\bx_0$ to the origin.
Then for small $t > t_{\rm slide}$, the above results imply
\begin{equation}
p_\ee(\bx,t;\bx_0) \to (1 - \varrho_{\rm no}) \delta(\bx - \bx^-(t - t_{\rm slide}))
+ \varrho_{\rm no} \delta(\bx - \bx^+(t - t_{\rm slide})) \;,
\end{equation}
as $\ee \to 0$.

\section{Conclusions}
\label{sec:CONC}
\setcounter{equation}{0}

The two-fold of Fig.~\ref{fig:schemPlanarTwoFold} has two important properties:
(i) many orbits reach the two-fold,
and, (ii) forward evolution from the two-fold is ambiguous.
To clarify the second property,
the two-fold has uncountably many distinct forward orbits in the sense of Filippov.
Specifically there are two tangent trajectories, $\bx^+(t)$ and $\bx^-(t)$,
and there are orbits that travel along the positive $y$-axis for some time,
then enter and evolve in either the left or right half-plane.
Local dynamics are described by the two-dimensional ODE system, (\ref{eq:ode}),
which has two parameters, $\cA$ and $\cB$,
that relate to the strength of the vector field for $x < 0$, relative to that for $x > 0$,
in the $x$ and $y$ directions respectively.

In this paper three different $O(\ee)$ perturbations to (\ref{eq:ode}) have been studied: hysteresis, time-delay and noise.
Each perturbation removes the ambiguity of forward evolution
and in each case forward evolution close to the two-fold has been analysed.
This microscopic analysis of the two-fold uncovers a macroscopic view of the dynamics.
For small $\ee > 0$, most orbits (or sample paths in the case of noise),
follow close to either $\bx^+(t)$ or $\bx^-(t)$.
In the limit, $\ee \to 0$, almost all orbits or sample paths follow exactly one of the these trajectories.
Thus the possibility of evolution along the positive $y$-axis for a nonzero length of time may be excluded.
This observation mirrors theoretical results for
stochastically perturbed continuous, non-differentiable systems in the zero-noise limit \cite{Ve83,BaBa82}.

With a hysteretic band of width $2 \ee$ the system is described by (\ref{eq:hysteresis}),
with a time-delay of $t=\ee$ the system is described by (\ref{eq:timedelay}),
and with $O(\ee)$ noise the system is described by (\ref{eq:sde}).
It is interesting that the three different types of perturbation exhibit three different scaling laws.
Specifically, for hysteresis, $x = O(\ee)$, for orbits approaching the two-fold about the negative $y$-axis.
For time-delay, $x = O(\ee^2)$, near the two-fold.
Finally for noise, $x = O \left( \ee^{\frac{4}{3}} \right)$, near the two-fold.

For perturbation by hysteresis and time-delay,
$\varrho_{\rm hy}$ and $\varrho_{\rm td}$, respectively,
represent the fraction of orbits that head right
(that is, follow close to $\bx^+(t)$ after passing by the two-fold)
in the limit $\ee \to 0$.
Similarly, for the stochastic perturbation, $\varrho_{\rm no}$ represents the limiting probability that sample paths head right.
As a function of $\cA$ and $\cB$,
an explicit expression was only obtained for $\varrho_{\rm hy}$, specifically (\ref{eq:varrhohysteresis2}).
For $\varrho_{\rm td}$ and $\varrho_{\rm no}$ efficient algorithms were established to evaluate these values numerically.
Lemma \ref{le:specialValues} gives special values of $\varrho_{\rm no}$.
Plots of each $\varrho(\cA,\cB)$ are shown in Figs.~\ref{fig:ABhysteresis}, \ref{fig:ABtimedelay} and \ref{fig:ABnoise}.
Each plot has the same scale and viewpoint so that they may be compared fairly.
For example, we notice that as $\cA \to 0$, for hysteresis and time-delay $\varrho \to 0$,
but for noise $\varrho \to \frac{1}{2}$.
Also $\varrho_{\rm no}(\cA,\cB)$ appears relatively flat, but, by (\ref{eq:sV1}), may take value in $(0,1)$.

The primary goal of this work has been to resolve the ambiguity of forward evolution from a two-fold.
To conclude, we may interpret forward evolution of the unperturbed system (\ref{eq:ode})
from the two-fold as probabilistic,
with evolution along $\bx^+(t)$ with probability, $\varrho(\cA,\cB)$,
and along $\bx^-(t)$ with probability, $1 - \varrho(\cA,\cB)$.
However, the value of $\varrho$ depends upon the type of perturbation under consideration.
This observation indicates that such a regularisation of two-folds is
potentially futile in the absence of additional information about the system.
As a mathematical model the system represents only an approximation to real-world behaviour.
For a given application, physical reasons may justify
the use of one of the three types of perturbation over the other two.
It remains to apply a probabilistic representation for forward evolution from a two-fold to a physical system.
Certainly by discounting Filippov solutions that traverse part of the repelling sliding region,
the two-piece probabilistic solution proposed here is more sophisticated
than permitting all possible Filippov solutions.
It remains to extend the ideas described here to two-folds in three or more dimensions
that arise generically in discontinuous systems of ODEs.

\appendix
\section{Proof of Lemma \ref{le:qI}}
\label{sec:HYSTPROOF}
\setcounter{equation}{0}

Here we write (\ref{eq:ykplus2}) and (\ref{eq:yk}) as
\begin{eqnarray}
y_k &=& -\sqrt{y_{k-2}^2 + 4 \left( 1 + \frac{\cB}{\cA} \right) \ee} + G(y_{k-2};\ee) \ee \;, {\rm ~for~} k \ge 3 \;, \label{eq:ykplus22} \\
y_k &=& \left\{ \begin{array}{lc}
-\sqrt{2} \sqrt{k + \frac{(k-1) \cB}{\cA}} \sqrt{\ee} + \zeta_k(\ee) \ee \;, & {\rm ~for~} k = 1,3,5,\ldots \\
-\sqrt{2} \sqrt{k-1 + \frac{k \cB}{\cA}} \sqrt{\ee} + \zeta_k(\ee) \ee \;, & {\rm ~for~} k = 2,4,6,\ldots
\end{array} \right.
\;. \label{eq:yk2}
\end{eqnarray}
There are two main parts to the proof.
First it is shown that there exists $M_1 \in \mathbb{R}$ such that
\begin{equation}
|\zeta_k(\ee)| \le k M_1 \;, \quad \forall k \le \ee^{-\theta} \;,
\label{eq:zetaBound}
\end{equation}
given $0 < \theta < 1$.
Second it is shown that
\begin{equation}
|q_{[y_{k+4},y_{k+2}]} - q_{[y_{k+2},y_k]}| \le M_2 \ee^{1+\theta} \;, \quad \forall k \ge \ee^{-\theta} \;,
\label{eq:qkplus2Bound}
\end{equation}
for some $M_2 \in \mathbb{R}$.
Lemma \ref{le:qI} is then a consequence of a few additional observations.
The key point is that the error in (\ref{eq:qkplus2Bound}) is of higher order than $\ee$
which enables us to usefully evaluate $q_{[y_k,y_{k-2}]}$ at $k = O(\ee^{-1})$.

We restrict our attention to $y \in [y^*,0]$, $\ee \in [0,\ee^*]$, where $y^*$ and $\ee^*$ are suitably small.
For this reason we may assume that the error in (\ref{eq:ykplus22}), $G$, is a smooth bounded function and let
$M_1 = 2 \,{\rm max}(|G|) + {\rm max}(|\zeta_1|) + \frac{1}{2} {\rm max}(|\zeta_2|)$.
Trivially (\ref{eq:zetaBound}) holds for $k = 1$ and $k = 2$.
Mathematical induction is required to verify (\ref{eq:zetaBound}) for all $k \le \ee^{-\theta}$.
Combining (\ref{eq:ykplus22}) and (\ref{eq:yk2}) leads to
\begin{eqnarray}
\zeta_{k+2} &=& G(y_k) + \frac{\sqrt{2}}{\sqrt{\ee}} \sqrt{k+2 + \frac{(k+1)\cB}{\cA}} \nonumber \\
&&-~\frac{\sqrt{2}}{\sqrt{\ee}}
\sqrt{ k+2 + \frac{(k+1)\cB}{\cA} - \sqrt{2} \sqrt{k + \frac{(k-1)\cB}{\cA}} \zeta_k \sqrt{\ee}
+ \frac{1}{2} \zeta_k^2 \ee} \;,
\label{eq:zetakplus2}
\end{eqnarray}
for odd values of $k$, and a similar formula for even $k$.
By appropriately bounding the large square-root term in (\ref{eq:zetakplus2})\removableFootnote{
To achieve this we use the following lemma which we state without proof:
\begin{lemma}
For any $a,b,c \in \mathbb{R}$, satisfying $a,c > 0$, $b^2 < 4 a c$, and for any $|x| < \frac{a}{|b|}$,
\begin{equation}
\left| \sqrt{a + bx + cx^2} - \sqrt{a} - \frac{b}{2 \sqrt{a}} x \right| \le \frac{4 a c - b^2}{4 a^{\frac{3}{2}}} x^2 \;.
\end{equation}
\label{le:squareRootBound}
\end{lemma}
We are able to apply Lemma \ref{le:squareRootBound} with:
$a = k+2 + \frac{(k+1)\cB}{\cA}$,
$b = -\sqrt{2} \sqrt{k + \frac{(k-1)\cB}{\cA}}$,
$c = \frac{1}{2}$,
and $x = \zeta_k \sqrt{\ee}$,
because, by assuming (\ref{eq:eeCondition}), we have
$|x| = |\zeta_k| \sqrt{\ee}
\le k M_1 \sqrt{\ee}
\le \sqrt{k} M_1 \ee^{\frac{1}{2}(1-\theta)}
< \frac{\sqrt{k}}{\sqrt{2}}
< \frac{\sqrt{k}}{\sqrt{2}} + \frac{\sqrt{2} \left( 1 + \frac{\cB}{\cA} \right)}{\sqrt{k + \frac{(k-1)\cB}{\cA}}}
= \frac{a}{|b|}$.
This leads to (\ref{eq:zetakplus2bound}).
},
which requires assuming $\ee^*$ is sufficiently small that\removableFootnote{
Really the $\left( 1 + \frac{\cB}{\cA} \right)$ piece is only required in the next step.
}
\begin{equation}
\ee^{\frac{1}{2}(1-\theta)} \le \frac{1}{2 \sqrt{2} \left( 1 + \frac{\cB}{\cA} \right) M_1} \;,
\label{eq:eeCondition}
\end{equation}
it may be shown that
\begin{equation}
|\zeta_{k+2} - \zeta_k| \le |G(y_k)| +
|\zeta_k| \left( 1 - \frac{\sqrt{k+\frac{(k-1)\cB}{\cA}}}{\sqrt{k+2+\frac{(k+1)\cB}{\cA}}} \right)
+ \frac{\sqrt{2} \left( 1 + \frac{\cB}{\cA} \right) \zeta_k^2 \sqrt{\ee}}{k^{\frac{3}{2}}} \;.
\label{eq:zetakplus2bound}
\end{equation}
By further using (\ref{eq:eeCondition}),
noting $\frac{1}{k} > 1 - \frac{\sqrt{k+\frac{(k-1)\cB}{\cA}}}{\sqrt{k+2+\frac{(k+1)\cB}{\cA}}}$,
and applying the induction hypothesis\removableFootnote{
$\frac{\sqrt{2} \left( 1 + \frac{\cB}{\cA} \right) \zeta_k^2 \sqrt{\ee}}{k^{\frac{3}{2}}}
\le \sqrt{2} \left( 1 + \frac{\cB}{\cA} \right) \sqrt{k} M_1^2 \sqrt{\ee}
\le \sqrt{2} \left( 1 + \frac{\cB}{\cA} \right) M_1^2 \ee^{\frac{1}{2}(1-\theta)}
\le M_1$
},
we obtain
\begin{equation}
|\zeta_{k+2} - \zeta_k| \le \frac{1}{2} M_1 + M_1 + \frac{1}{2} M_1 = 2 M_1 \;.
\label{eq:zetakplus2bound2}
\end{equation}
Consequently, $|\zeta_{k+2}| \le (k+2) M_1$, for odd values of $k$.
The result for even $k$ may be shown in the same manner.
This completes the inductive demonstration of (\ref{eq:zetaBound}).

Since (\ref{eq:qI}) involves only odd values of $k$,
for remainder of the proof it is assumed that $k$ is odd
(even values of $k$ may be treated similarly).
Substituting (\ref{eq:yk2}) into (\ref{eq:qk}) yields
\begin{equation}
q_{[y_k,y_{k-2}]} =
\frac{\cA + O \left( \frac{1}{k} \right)
+ \frac{\cA \sqrt{k \ee}}{\sqrt{2}} \sqrt{1 + \frac{\cB}{\cA}} \left( \zeta_{k-1} - \zeta_k \right)}
{\cA+\cB + O \left( \frac{1}{k} \right)
+ \frac{\cA \sqrt{k \ee}}{\sqrt{2}} \sqrt{1 + \frac{\cB}{\cA}} \left( \zeta_{k-2} - \zeta_k \right)} \;.
\label{eq:qk3}
\end{equation}
By applying (\ref{eq:zetaBound}), an evaluation of (\ref{eq:qk3}) at $k = O(\ee^{-\theta})$ gives
\begin{equation}
q_{[y_k,y_{k-2}]} =
\frac{\cA}{\cA+\cB} + O(\ee^\theta) + O \left( \ee^{\frac{1}{2} (1-3\theta)} \right) \;.
\label{eq:qktheta}
\end{equation}
To minimise the magnitude of the error we take $\theta = \frac{1}{5}$.
It is likely that a stronger bound on the error is achievable,
but (\ref{eq:qktheta}) is sufficient for this proof.
The primary difficulty is bounding $\zeta_{k-1} - \zeta_k$, in (\ref{eq:qk3}), for which
(\ref{eq:zetaBound}) was used crudely to obtain $|\zeta_{k-1} - \zeta_k| \le 2 k M_1$.

Next, (\ref{eq:qkplus2Bound}) is proved.
By (\ref{eq:yk2}), for all $k \ge \ee^{-\theta}$ we have $|y_k| > \ee^{\frac{1}{2}(1-\theta)}$.
Thus\removableFootnote{
We have
\begin{equation}
|y_{k+2}-y_k| \le |G(y_k)| \ee + \left| y_k + \sqrt{y_k^2 + 4 \left( 1 + \frac{\cA}{\cB} \right) \ee} \right| \;.
\end{equation}
Since the last term in this expression is an increasing function of $y_k$,
we substitute $y_k = \ee^{\frac{1}{2}(1-\theta)}$,
and apply simple bound to the square root to arrive at
\begin{equation}
|y_{k+2}-y_k| \le 4 \left( 1 + \frac{\cB}{\cA} \right) \ee^{\frac{1}{2}(1+\theta)} \;,
\end{equation}
for sufficiently small $\ee$.
}
by (\ref{eq:ykplus22}), we can write
\begin{equation}
y_{k+2} = y_k - H(y_k;\ee) \ee^{\frac{1}{2}(1+\theta)} \;,
\label{eq:ykplus23}
\end{equation}
where $H$ is non-negative, bounded and has bounded derivatives.
By (\ref{eq:ykplus23}) and (\ref{eq:qk}),
\begin{equation}
y_{k+1} = y_k - (1-q_{[y_{k+2},y_k]}) H(y_k) \ee^{\frac{1}{2}(1+\theta)} \;.
\label{eq:ykplus1}
\end{equation}
Further applications of (\ref{eq:ykplus23}) and (\ref{eq:qk}) with (\ref{eq:ykplus1}) lead to
\begin{equation}
q_{[y_{k+4},y_{k+2}]} =
\frac{H(y_{k+2}) - H(y_{k+1}) + q_{[y_{k+2},y_k]} H(y_k)}{H(y_{k+2})} \;.
\label{eq:qk4}
\end{equation}
By Taylor expanding the various instances of $H$ in (\ref{eq:qk4}) about $y_k$,
we find that all $O(1)$ and $O(\ee^{\frac{1}{2}(1+\theta)})$ terms cancel producing,
$q_{[y_{k+4},y_{k+2}]} - q_{[y_{k+2},y_k]} = O(\ee^{1+\theta})$.
Moreover, assuming the vector field, (\ref{eq:hysteresis}), is at least piecewise-$C^1$,
there exists a sufficiently small range of $y$ and $\ee$ values over which $H''(y_k)$ is bounded,
which establishes (\ref{eq:qkplus2Bound}).

Finally, for $\ee$-independent values, $y_{\rm min}$ and $y_{\rm max}$, $k_{\rm min}$ and $k_{\rm max}$ are $O(\ee^{-1})$.
Thus, by (\ref{eq:qkplus2Bound}) and (\ref{eq:qktheta}) with $\theta = \frac{1}{5}$,
\begin{equation}
q_{[y_{k_{\rm min}+2},y_{k_{\rm min}}]} = \frac{\cA}{\cA+\cB} + O(\ee^{\frac{1}{5}}) \;,
\end{equation}
and for any odd $k$ with, $k_{\rm min}+2 \le k \le k_{\rm max}$,
the difference, $q_{[y_k,y_{k-2}]} - q_{[y_{k_{\rm min}+2},y_{k_{\rm min}}]}$,
is $O(\ee^{\frac{1}{5}})$.
Therefore, by (\ref{eq:qIbound}),
\begin{equation}
q_{[y_{\rm min},y_{\rm max}]} = \frac{\cA}{\cA+\cB} + O(\ee^{\frac{1}{5}}) \;,
\label{eq:qI3}
\end{equation}
which completes the proof.
(Additional arguments to maximise the exponent in the error term of (\ref{eq:qI3}) are beyond the scope of this paper.)

\section{Derivation of (\ref{eq:sV1})}
\label{sec:SV1}
\setcounter{equation}{0}

If $\cB = \cA^2$, then by stretching the negative $u_0$-axis by a factor, $\frac{1}{\cA}$,
we find that (\ref{eq:Qpde2}) behaves similarly for negative and positive values of $u_0$.
To exploit this symmetry, let
\begin{equation}
Z(u_0,r) = Q(u_0,r) + \cA Q \left( -\frac{u_0}{\cA},r \right) \;.
\label{eq:Z}
\end{equation}
From (\ref{eq:Qpde2})-(\ref{eq:bcs2}),
for $u_0 > 0$ and $r \in \mathbb{R}$, $Z$ satisfies the boundary value problem
\begin{equation}
\frac{\partial Z}{\partial r} =
-r \frac{\partial Z}{\partial u_0} +
\frac{1}{2} \frac{\partial^2 Z}{\partial u_0^2} \;,
\end{equation}
\begin{equation}
\lim_{r \to -\infty} Z(u_0,r) = 1 \;, \qquad
\frac{\partial Z}{\partial u_0}(0,r) = 0 \;, \qquad
\lim_{u_0 \to \infty} Z(u_0,r) = 1 \;.
\end{equation}
The unique solution to this problem is $Z(u_0,r) \equiv 1$.
By (\ref{eq:Z}), for all $r$, $Q(0,r) = \frac{1}{1+\cA} Z(0,r) = \frac{1}{1+\cA}$.
Thus, in particular, $\varrho_{\rm no} = \lim_{r \to \infty} Q(0,r) = \frac{1}{1+\cA}$,
as required.

\section{Derivation of (\ref{eq:sV2})}
\label{sec:SV2}
\setcounter{equation}{0}

If $\cB = 1$, then (\ref{eq:sdeReduced}) has two special properties that provide simplification.
First, $v(s)$ increases at a unit rate at all times.
Assuming, for simplicity, $v(0) = 0$, we therefore have
\begin{equation}
du(s) = \left\{ \begin{array}{lc}
-\cA s \;, & u < 0 \\
s \;, & u > 0 
\end{array} \right\} \,ds + dW(s) \;,
\label{eq:sdeB1}
\end{equation}
and $\varrho_{\rm no} = \lim_{s \to \infty} \mathbb{P} [ u(s) > 0 ]$\removableFootnote{
This is intuitive but non-trivial.
The claim is that $Q = \lim_{s \to \infty} \mathbb{P} [ u(s) > 0 ]$.
}
Second, we may write (\ref{eq:sdeB1}) as
\begin{equation}
du(s) = \phi(u,s) \,ds + dW(s) \;, 
\label{eq:sdeB1gen}
\end{equation}
where
\begin{equation}
\phi(u,s) = \left\{ \begin{array}{lc}
-\cA s \;, & u < 0 \\
s \;, & u > 0
\end{array} \right. \;.
\label{eq:specificphi}
\end{equation}
The key observation is that the drift in (\ref{eq:specificphi}) is odd with respect to time, that is,
\begin{equation}
\phi(u,-s) = -\phi(u,s) \;,
\label{eq:phiSymmetry}
\end{equation}
for all $u$ and $s$.
For this reason, roughly speaking, drift over negative times is cancelled by drift over positive times
and consequently the probability that $u(s)$ is positive approaches $\frac{1}{2}$ as $s \to \infty$.
Indeed (\ref{eq:sV2}) is a consequence of the following general result:

\begin{lemma}
Suppose $\phi$ is a bounded, measurable function satisfying (\ref{eq:phiSymmetry})
with possibly finitely many values of $u$ at which it is discontinuous.
For any $T > 0$, for the stochastic process (\ref{eq:sdeB1gen}),
$\mathbb{P} \left[ u(T) > 0 ~\big|~ u(-T) = 0 \right] = \frac{1}{2}$.
\label{le:generalSymmetry}
\end{lemma}

The desired result, (\ref{eq:sV2}), follows by applying Lemma \ref{le:generalSymmetry} to (\ref{eq:specificphi})
and taking $T \to \infty$.
Below, Lemma \ref{le:generalSymmetry} is proved for smooth $\phi$.
For brevity a generalisation permitting discontinuities in $\phi$ is omitted.
This may be achieved by simply adding consistency conditions
for the transitional probability density function of (\ref{eq:sdeB1gen}) at points of discontinuity.

For (\ref{eq:sdeB1gen}), the probability density function of $u(s)$, call it $p(u,s)$,
is the unique bounded solution to
\begin{equation}
\frac{\partial p}{\partial s} =
-\frac{\partial (\phi p)}{\partial u} +
\frac{1}{2} \frac{\partial^2 p}{\partial u^2} \;, \qquad
p(u,-T) = \delta(u) \;.
\end{equation}
Here we let
\begin{equation}
Q(u_0,s_0) \equiv \mathbb{P} \left[
u(T) > 0 ~\big|~ u(s_0) = u_0 \right] \;.
\end{equation}
By the Feynman-Kac formula \cite{Ok03,KaSh91}\removableFootnote{
{\bf The Feynman-Kac formula}:
With certain smoothness and growth constraints on functions $a$, $b$, $f$ and $g$,
for the stochastic differential equation
\begin{equation}
du(s) = a(u(s),s) \,ds + b(u(s),s) \,dW(s) \;,
\end{equation}
if
\begin{equation}
Q(u_0,s_0) = \mathbb{E}_{u_0} \left( f(u(s)) \,{\rm e}^{-\int_0^s g(u(\hat{s})) \,d\hat{s}} \right) \;,
\end{equation}
then $Q(u_0,s_0)$ is the unique solution to
\begin{equation}
Q_{s_0} = -a(u_0,s_0) Q_{u_0} - \frac{b(u_0,s_0)^2}{2} Q_{u_0 u_0} + g(u_0) Q \;,
\end{equation}
with final condition, $Q(u_0,T) = f(u_0)$,
on the domain $\mathbb{R} \times (-\infty,T)$.
For our problem, $a = \phi$, $b = 1$, $f = \chi_{[0,\infty)}$, and $g = 0$.
},
$Q$ is the unique, bounded, $C^1$ solution to
\begin{equation}
\frac{\partial Q}{\partial s_0} =
-\phi(u_0,s_0) \frac{\partial Q}{\partial u_0}
- \frac{1}{2} \frac{\partial^2 Q}{\partial u_0^2} \;, \qquad
Q(u_0,T) = \chi_{[0,\infty)}(u_0) \;.
\end{equation}
By the symmetry property, (\ref{eq:phiSymmetry}),
the solutions to these two partial differential equations are related by
$\frac{\partial Q}{\partial u_0}(u_0,s_0) = p(u_0,-s_0)$,
for all $u_0$ and $s_0$.
Now define
\begin{equation}
h(s_0) \equiv \int_{-\infty}^\infty
p(u,-s_0) Q(u,s_0) \,du \;.
\end{equation}
A quick manipulation shows that this function takes the value $\frac{1}{2}$
for any value of $s_0$:
\begin{equation}
h(s_0) = \int_{-\infty}^\infty
\frac{\partial Q}{\partial u_0}(u,s_0) Q(u,s_0) \,du
= \frac{1}{2} \int_{-\infty}^\infty
\frac{\partial}{\partial u} \left( Q(u,s_0) \right)^2 \,du
= \frac{1}{2} \;,
\end{equation}
since for any fixed $s_0$,
$Q(u_0,s_0) \to 0$ as $u_0 \to -\infty$ and
$Q(u_0,s_0) \to 1$ as $u_0 \to \infty$.
Also, by the strong Markov property, for any $s_0$,
\begin{equation}
\mathbb{P} \left[ u(T) > 0 ~\big|~ u(-T) = 0 \right] =
\int_{-\infty}^\infty p(u,s_0) Q(u,s_0) \,du \;.
\end{equation}
By then setting $s_0 = 0$ we deduce 
\begin{equation}
\mathbb{P} \left[ u(T) > 0 ~\big|~ u(-T) = 0 \right] = h(0) = \frac{1}{2} \;,
\end{equation}
as required.

\section{Derivation of (\ref{eq:sV3})}
\label{sec:SV3}
\setcounter{equation}{0}

It easier to treat the limit $\cA \to \infty$.
Here we show that
\begin{equation}
\lim_{\cA \to \infty} \varrho_{\rm no}(\cA,\cB) = \frac{1}{2} \;,
\label{eq:sV3alt}
\end{equation}
for all $\cB > 0$.
In view of (\ref{eq:varrhoSym}), this is equivalent to (\ref{eq:sV3}).

For (\ref{eq:sdeReduced}), the limit, $\cA \to \infty$, represents infinitely strong drift for $u < 0$.
Therefore a sample path from a point $(u,v)$ with $u < 0$ and $v < 0$ is immediately sent to $(0,v)$, with probability $1$.
If instead $u < 0$ and $v > 0$, then the $u$-value of the solution goes to $-\infty$, with probability $1$.
Consequently, it suffices to analyse dynamics for $u > 0$
and impose a reflecting boundary condition on the negative $v$-axis
and an absorbing boundary condition on the positive $v$-axis.
Specifically, we may assume $v(s) \equiv s$, and that $u(s)$ is governed by
\begin{equation}
du(s) = s \,ds + dW(s) \;.
\label{eq:sdeAinf}
\end{equation}
The probability density function for $u(s)$, call it $p(u,s)$, therefore satisfies
\begin{equation}
\frac{\partial p}{\partial s} =
-s \frac{\partial p}{\partial u} +
\frac{1}{2} \frac{\partial^2 p}{\partial u^2} \;,
\label{eq:fpeR}
\end{equation}
\begin{equation}
\frac{\partial p}{\partial u}(0,s) = 2 s p(0,s) \;, {\rm ~for~} s < 0 \;, \qquad
p(0,s) = 0 \;, {\rm ~for~} s > 0 \;,
\label{eq:bcsR1}
\end{equation}
\begin{equation}
\lim_{s \to -\infty} p(u,s) = \delta(u) \;, \qquad
\lim_{u \to \infty} p(u,s) = 0 \;.
\label{eq:bcsR2}
\end{equation}

Note, in particular, the behaviour of (\ref{eq:sdeReduced}) in the limit, $\cA \to \infty$, is independent of $\cB$.
Since, by (\ref{eq:sV2}), $\lim_{\cA \to \infty} \varrho_{\rm no}(\cA,1) = \frac{1}{2}$,
we therefore have (\ref{eq:sV3alt}) for all $\cB$, as required.
However, to add to our understanding of the system for large (and small) values of $\cA$,
here we derive an explicit expression of $p(u,s)$ from (\ref{eq:fpeR})-(\ref{eq:bcsR2}),
and use it to derive (\ref{eq:sV3alt}) directly.

The problem, (\ref{eq:fpeR})-(\ref{eq:bcsR2}),
with instead a reflecting boundary condition at $u = 0$ for all values of $s$,
was solved by Knessl in \cite{Kn00}.
Thus for $s \le 0$, $p(u,s)$ is Knessl's solution:
\begin{equation}
p(u,s) = 2^{\frac{2}{3}} \,{\rm e}^{-\frac{s^3}{6} + s u} K(u,s) \;,
\label{eq:pLower}
\end{equation}
where $K$ may be expressed as an inverse Laplace transform
\begin{equation}
K(u,s) = \frac{1}{2 \pi {\rm i}} \int_{\rm Br}
{\rm e}^{\alpha s}
\frac{{\rm Ai} \left( 2^{\frac{1}{3}}(u + \alpha) \right)}
{{\rm Ai}^2 \left( 2^{\frac{1}{3}} \alpha \right)} \,d\alpha \;.
\label{eq:K}
\end{equation}
In (\ref{eq:K}), ${\rm Ai}(z)$ is the Airy function of the first kind, and ${\rm Br}$ denotes a suitable Bromwich contour in the complex plane.
Throughout this exposition Bromwich contours may be assumed to be vertical lines at a non-negative real value.
The following expression of $p$ for $s \ge 0$ was obtained by methods similar to those in \cite{Kn00}
(for brevity the details are omitted)\removableFootnote{
If we write
\begin{equation}
p(u,s) = {\rm e}^{-\frac{s^3}{6} + su} \psi(u,s) \;,
\end{equation}
then $\psi$ is specified by the boundary value problem
\begin{equation}
\psi_s = -u \psi + \frac{1}{2} \psi_{uu} \;,
\end{equation}
\begin{equation}
\psi(0,s) = 0 \;, \qquad
\lim_{u \to \infty} \psi(u,s) = 0 \;, \qquad
\psi(u,0) = 2^{\frac{2}{3}} K(u,0) \;.
\end{equation}
Let
\begin{equation}
\Psi(u,\beta) = \int_0^\infty {\rm e}^{-\beta s} \psi(u,s) \,ds \;,
\end{equation}
then
\begin{equation}
\frac{1}{2} \Psi_{uu} - (u+\beta) \Psi = -2^{\frac{2}{3}} K(u,0) \;,
\end{equation}
\begin{equation}
\Psi(0,\beta) = 0 \;, \qquad
\lim_{u \to \infty} \Psi(u,\beta) = 0 \;.
\end{equation}
Thus we may write
\begin{equation}
\Psi(u,\beta) = \int_0^\infty -2^{\frac{2}{3}} K(\xi,0) \Phi(u,\beta,\xi) \,d\xi \;,
\end{equation}
where
\begin{equation}
\frac{1}{2} \Phi_{uu} - (u+\beta) \Phi = \delta(u-\xi) \;,
\end{equation}
\begin{equation}
\Phi(0,\beta,\xi) = 0 \;, \qquad
\lim_{u \to \infty} \Phi(u,\beta,\xi) = 0 \;.
\end{equation}
Standard ODE methods may then be used to derive (\ref{eq:Phi}).
}:
\begin{equation}
p(u,s) = 2^{\frac{2}{3}} {\rm e}^{-\frac{s^3}{6} + s u}
\frac{1}{2 \pi {\rm i}} \int_{\rm Br} {\rm e}^{\beta s}
\int_0^\infty -K(\xi,0) \Phi(u,\beta,\xi) \,d\xi \,d\beta \;,
\label{eq:pUpper}
\end{equation}
where
\begin{equation}
\Phi(u,\beta,\xi) = \left\{ \begin{array}{lc}
\frac{2^{\frac{2}{3}} \pi {\rm Ai} \left( 2^{\frac{1}{3}} (\xi+\beta) \right)}
{{\rm Ai} \left( 2^{\frac{1}{3}} \beta \right)}
\left( {\rm Ai} \left( 2^{\frac{1}{3}} (u+\beta) \right) {\rm Bi} \left( 2^{\frac{1}{3}} \beta \right)
+ {\rm Ai} \left( 2^{\frac{1}{3}} \beta \right) {\rm Bi} \left( 2^{\frac{1}{3}} (u+\beta) \right) \right) \;, & u < \xi \\
\frac{2^{\frac{2}{3}} \pi {\rm Ai} \left( 2^{\frac{1}{3}} (u+\beta) \right)}
{{\rm Ai} \left( 2^{\frac{1}{3}} \beta \right)}
\left( {\rm Ai} \left( 2^{\frac{1}{3}} (\xi+\beta) \right) {\rm Bi} \left( 2^{\frac{1}{3}} \beta \right)
+ {\rm Ai} \left( 2^{\frac{1}{3}} \beta \right) {\rm Bi} \left( 2^{\frac{1}{3}} (\xi+\beta) \right) \right) \;, & u > \xi
\end{array} \right. \;,
\label{eq:Phi}
\end{equation}
and ${\rm Bi}(z)$ is the Airy function of the second kind.
We may use (\ref{eq:pUpper}) to derive (\ref{eq:sV3alt}) because\removableFootnote{
An alternative is to use,
$\lim_{\cA \to \infty} \varrho_{\rm no}(\cA,\cB) = 1 - \frac{1}{2} \int_0^\infty p_u(0,s) \,ds$,
where we have
\begin{equation}
p_u(0,s) = 2^{\frac{5}{3}} {\rm e}^{-\frac{s^3}{6}}
\int_0^\infty K(\xi,0)
\frac{1}{2 \pi {\rm i}}
\int_{\rm Br} 
{\rm e}^{\beta s}
\frac{{\rm Ai} \left( 2^{\frac{1}{3}}(\xi + \beta) \right)}
{{\rm Ai} \left( 2^{\frac{1}{3}} \beta \right)} \,d\alpha \,d\xi \;,
\end{equation}
but we have been unable to gain use from this.
}
\begin{equation}
\lim_{\cA \to \infty} \varrho_{\rm no}(\cA,\cB) =
\lim_{s \to \infty} \int_0^\infty p(u,s) \,du \;.
\label{eq:sV3alt2}
\end{equation}

For large $s$, due to the nature of Brownian motion
a solution to (\ref{eq:sdeAinf}) that has not been absorbed at $u=0$
will very likely lie an $O \left( s^{\frac{1}{2}} \right)$ distance from
the value $\frac{1}{2} s^2$ (the solution to (\ref{eq:sdeAinf}) in the absence of the noise term).
For this reason we write
\begin{equation}
u = \frac{1}{2} s^2 + \Delta s^{\frac{1}{2}} \;,
\label{eq:Delta}
\end{equation}
and take $u,s \to \infty$ with $\Delta = O(1)$.

As $u \to \infty$, only the piece of (\ref{eq:Phi}) with $u > \xi$ is important.
From this we utilise various basic properties of Airy functions to obtain, for large $u$\removableFootnote{
Here we use the Wronskian identity,
\begin{equation}
{\rm Ai}(z) {\rm Bi}'(z) - {\rm Ai}'(z) {\rm Bi}(z) = \frac{1}{\pi} \;,
\end{equation}
and 
\begin{eqnarray}
\int_0^\infty {\rm Ai} \left( 2^{\frac{1}{3}} (\xi+\alpha) \right)
{\rm Ai} \left( 2^{\frac{1}{3}} (\xi+\beta) \right) \,d\xi &=&
\frac{{\rm Ai}' \left( 2^{\frac{1}{3}} \alpha \right) {\rm Ai} \left( 2^{\frac{1}{3}} \beta \right) -
{\rm Ai} \left( 2^{\frac{1}{3}} \alpha \right) {\rm Ai}' \left( 2^{\frac{1}{3}} \beta \right)}
{2^{\frac{2}{3}} (\beta - \alpha)} \;, \\
\int_0^\infty {\rm Ai} \left( 2^{\frac{1}{3}} (\xi+\alpha) \right)
{\rm Bi} \left( 2^{\frac{1}{3}} (\xi+\beta) \right) \,d\xi &=&
\frac{{\rm Ai}' \left( 2^{\frac{1}{3}} \alpha \right) {\rm Bi} \left( 2^{\frac{1}{3}} \beta \right) -
{\rm Ai} \left( 2^{\frac{1}{3}} \alpha \right) {\rm Bi}' \left( 2^{\frac{1}{3}} \beta \right)}
{2^{\frac{2}{3}} (\beta - \alpha)} \;,
\end{eqnarray}
where the last identity requires (\ref{eq:Bromwichassumption}).
},
\begin{equation}
\int_0^\infty -K(\xi,0) \Phi(u,\beta,\xi) \,d\xi ~\sim~
{\rm Ai} \left( 2^{\frac{1}{3}} (u+\beta) \right) F(\beta) \;,
\label{eq:firstAsymptApprox}
\end{equation}
where
\begin{equation}
F(\beta) = \frac{1}{2 \pi {\rm i}} \int_{\rm Br}
\frac{1}{{\rm Ai} \left( 2^{\frac{1}{3}} \alpha \right) {\rm Ai} \left( 2^{\frac{1}{3}} \beta \right)
(\alpha - \beta)} \,d\alpha \;,
\label{eq:F}
\end{equation}
subject to the restriction
\begin{equation}
{\rm Re}(\alpha) > {\rm Re}(\beta) \;.
\label{eq:Bromwichassumption}
\end{equation}

Next, to (\ref{eq:pUpper}) with (\ref{eq:Delta}) and (\ref{eq:firstAsymptApprox}),
as in \cite{Kn00} we apply the formula,
${\rm Ai}(z) \sim \frac{1}{2 \pi^{\frac{1}{2}} z^{\frac{1}{4}}}
{\rm e}^{-\frac{2}{3} z^{\frac{3}{2}}}$,
valid for large $|z|$ bounded away from the negative real axis \cite{VaSo10},
which yields
\begin{equation}
p(u,s) \sim \frac{1}{\sqrt{2 \pi s}}
\,{\rm e}^{-\frac{\Delta^2}{2}}
\frac{2^{\frac{1}{3}}}{2 \pi {\rm i}} \int_{\rm Br} F(\beta) \,d\beta \;.
\end{equation}
Therefore as $s \to \infty$, $p$ is a scaled Gaussian and, by (\ref{eq:sV3alt2}) and (\ref{eq:F}),
\begin{equation}
\lim_{\cA \to \infty} \varrho_{\rm no}(\cA,\cB) =
\frac{2^{\frac{1}{3}}}{2 \pi {\rm i}} \int_{\rm Br} F(\beta) \,d\beta
= \frac{1}{(2 \pi {\rm i})^2} \int_{\rm Br} \int_{\rm Br}
\frac{1}{{\rm Ai}(a) {\rm Ai}(b) (a-b)} \,da \,db \;,
\label{eq:doubleInt}
\end{equation}
where we have substituted, $a = 2^{\frac{1}{3}} \alpha$,
and, $b = 2^{\frac{1}{3}} \beta$, and hence ${\rm Re}(a) > {\rm Re}(b)$.
By further writing, $a = \gamma + \frac{{\rm i} \tau}{2} + {\rm i} \omega$ (for any $\gamma > 0$),
and $b = -\frac{{\rm i} \tau}{2} + {\rm i} \omega$,
we produce
\begin{equation}
\lim_{\cA \to \infty} \varrho_{\rm no}(\cA,\cB) =
\frac{-{\rm i}}{4 \pi^2} \int_{-\infty}^\infty \int_{-\infty}^\infty
\frac{1}{{\rm Ai} \left( \gamma + \frac{{\rm i} \tau}{2} + {\rm i} \omega \right)
{\rm Ai} \left( -\frac{{\rm i} \tau}{2} + {\rm i} \omega \right)
(\tau - {\rm i} \gamma)} \,d\tau \,d\omega \;.
\label{eq:doubleInt2}
\end{equation}
Since (\ref{eq:doubleInt2}) is independent of the value of $\gamma$,
to evaluate it we may take the limit $\gamma \to 0$.
To achieve this we apply the following version of Plemelj-Sokhotski's formula \cite{AbFo03}:
$\lim_{\gamma \to 0^+} \int \frac{h(x)}{x - {\rm i} \gamma} \,dx
= {\rm i} \pi h(0) + {\rm C.P.} \int \frac{h(x)}{x} \,dx$,
(for integration over an interval containing zero),
which produces
\begin{equation}
\lim_{\cA \to \infty} \varrho_{\rm no}(\cA,\cB) =
\frac{1}{4 \pi} \int_{-\infty}^\infty \frac{1}{{\rm Ai}^2({\rm i}\omega)} \,d\omega =
-\frac{{\rm i}}{4} \,\frac{{\rm Bi}({\rm i}\omega)}{{\rm Ai}({\rm i}\omega)} \bigg|_{-\infty}^\infty =
\frac{1}{2} \;,
\end{equation}
as required.

\end{document}